\documentclass[a4paper,11pt]{article}
\usepackage{amsmath,amssymb, amsbsy, amsfonts, amsthm}
\usepackage{mathrsfs}
\usepackage{bigints}
\usepackage{color,psfrag}
\usepackage[dvips]{graphicx}
\usepackage{enumerate}
\usepackage{tikz}
\usetikzlibrary{decorations.markings}
\tikzset{
  set arrow inside/.code={\pgfqkeys{/tikz/arrow inside}{#1}},
  set arrow inside={end/.initial=>, opt/.initial=},
  /pgf/decoration/Mark/.style={
    mark/.expanded=at position #1 with
    {
      \noexpand\arrow[\pgfkeysvalueof{/tikz/arrow inside/opt}]{\pgfkeysvalueof{/tikz/arrow inside/end}}
    }
  },
  arrow inside/.style 2 args={
    set arrow inside={#1},
    postaction={
      decorate,decoration={
        markings,Mark/.list={#2}
      }
    }
  },
}
\tikzset{
  set arrow outside/.code={\pgfqkeys{/tikz/arrow outside}{#1}},
  set arrow outside={end/.initial=<, opt/.initial=},
  /pgf/decoration/Mark/.style={
    mark/.expanded=at position #1 with
    {
      \noexpand\arrow[\pgfkeysvalueof{/tikz/arrow outside/opt}]{\pgfkeysvalueof{/tikz/arrow outside/end}}
    }
  },
  arrow outside/.style 2 args={
    set arrow outside={#1},
    postaction={
      decorate,decoration={
        markings,Mark/.list={#2}
      }
    }
  },
}
\usepackage{url}
\usepackage{comment}
\usepackage{colortbl}

\usepackage{graphicx,color}
\graphicspath{{Figures/}}
\usepackage{tikz,pgfplots,pgf}
\usetikzlibrary{patterns}
\usepackage[font=small]{caption}
\usepackage{subcaption}

\newcommand{\eps}{\varepsilon}

\renewcommand{\geq }{\geqslant}

\renewcommand{\leq }{\leqslant}

\def\neweq#1{\begin{equation}\label{#1}}
\def\endeq{\end{equation}}

\newtheorem{theorem}{Theorem}
\newtheorem{proposition}[theorem]{Proposition}
\newtheorem{lemma}[theorem]{Lemma}
\newtheorem{corollary}[theorem]{Corollary}
\newtheorem{definition}[theorem]{Definition}

\newtheorem{remark}[theorem]{Remark}

\titlepage 
\title{Asymptotic study of critical wave fronts \\ for parameter-dependent Born-Infeld models: \\ physically predicted behaviors and new phenomena}
\author{Maurizio Garrione\footnote{The author acknowledges the financial support of the PRIN Project 201758MTR2 - \emph{Direct and inverse problems for partial differential equations: theoretical aspects and applications} (Italy) and of the Gruppo Nazionale per l'Analisi Matematica, la Probabilit\`a e le loro Applicazioni (GNAMPA) of the Istituto Nazionale di Alta Matematica (INdAM), Rome, Italy.}
}
\date{\small{Dipartimento di Matematica, Politecnico di Milano \\ 
Piazza Leonardo da Vinci 32, 20133 Milano \\
e-mail address: \tt{maurizio.garrione@polimi.it}}}

\begin{document}
\maketitle

\begin{abstract}
In this paper, we study some parameter-dependent reaction-diffusion models
governed by the \emph{Born-Infeld} (or \emph{Minkowski}) operator. In dependence on two parameters $a, b > 0$, related to the field strength and to the diffusivity, we investigate the limit critical speed for traveling fronts, together with the limit behavior of the associated critical profiles. As the two main results, on the one hand 
we rigorously show that for arbitrarily large electric fields the critical speed and the critical profiles converge to the ones of the linear diffusion problem, agreeing with the well-known physical prediction that, in this case, Born's law and Maxwell's law coincide. Such a result is accompanied by a counterpart for arbitrarily small electric fields and a complete analysis for vanishing/large diffusion.
On the other hand, we prove the onset of a new and unexpected phenomenon for the singular perturbation problem: the critical speed of propagation \emph{does not converge to zero} and, for KPP or combustion-type reactions, the limit front profile becomes sharp \emph{on one side only}. In fact, this latter coincides with the $C^1$-gluing of a piecewise linear profile having slope equal to $1$ when nonconstant, and of a regular inviscid profile propagating at the limit speed. The gluing point and the value of the limit speed are determined explicitly. The outcomes of the whole discussion, based on a careful analysis of the first-order reduction associated with the original equation, show substantial differences and peculiarities of the Born-Infeld operator with respect to linear or saturating diffusions, and are supported by several numerical experiments. \end{abstract}

\smallbreak

\section{Introduction}

\begin{figure}[ht!]
\includegraphics[width=\textwidth, height=5.1cm]{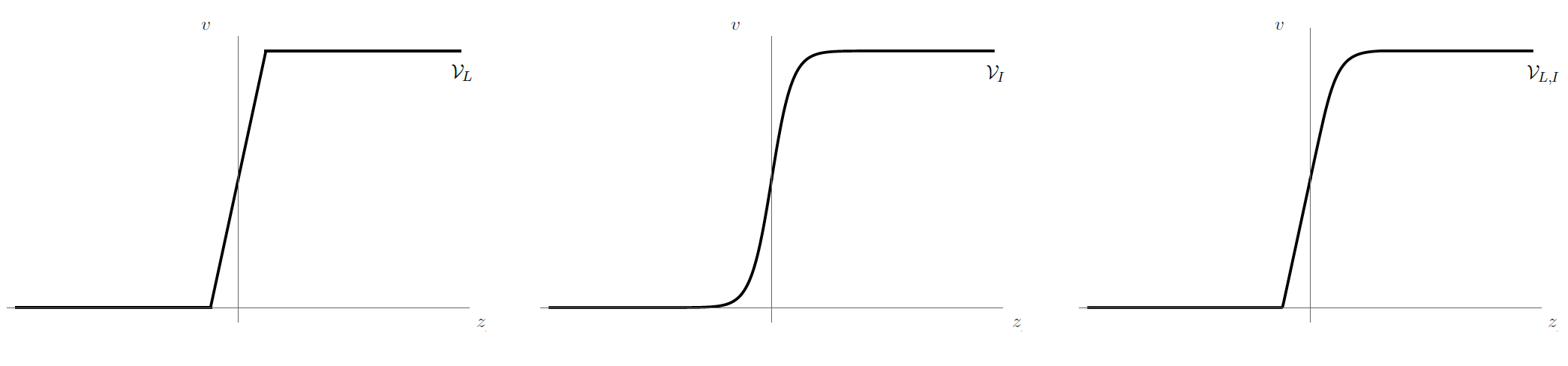}
\caption{The shapes of $\mathcal{V}_L$, $\mathcal{V}_I$ and $\mathcal{V}_{L, I}$}\label{introduttiva}
\end{figure}

In this paper, we discuss some features of traveling fronts for the one-dimensional nonlinear reaction-diffusion equation
\begin{equation}\label{MinkIntro}
u_t=\left(\frac{u_x}{\sqrt{a^2-b^2 u_x^2}}\right)_x  + f(u), \quad u=u(x, t), \; x \in \mathbb{R}, \, t \in \mathbb{R},
\end{equation}
where $a$ and $b$ are two positive parameters and $f(0)=0=f(1)$. Equation \eqref{MinkIntro} is 
inspired by the $n$-dimensional evolution equation
\begin{equation}\label{Mink}
u_t= \textnormal{div}\,\left(\frac{\nabla u}{\sqrt{1-\vert \nabla u \vert^2/\beta^2}}\right) + f(u), \quad u=u(x, t), \; x \in \mathbb{R}^n, \, t \in \mathbb{R},
\end{equation}
which naturally appears when modeling electromagnetic phenomena in presence of an upper bound ($\beta > 0$) on the value of the electric field. More precisely, one usually starts from the Born-Infeld Lagrangian $L=\beta^2 (1-\sqrt{1-(1/\beta^2)(E^2-B^2)})$ \cite{BorInf}; 
in the static case, assuming that the magnetic field $B$ is zero, by Maxwell equations the electric field is conservative, so that $E$ coincides with the gradient of an electrostatic potential $u$. If we admit the additional presence of a nonlinear potential density function $F(u)$ in the expression of the Lagrangian, 
the associated Euler-Lagrange equation for $u$ takes then the form 
\begin{equation}\label{lastazionaria}
\textnormal{div}\,\left(\frac{\nabla u}{\sqrt{1-\vert \nabla u \vert^2/\beta^2}}\right) + f(u) = 0 
\end{equation}
(where $f=F'$), namely is the stationary counterpart of \eqref{Mink}. Such an equation has deep connections with string theory, see for instance \cite{Gib, Yan};
the positive parameter $\beta$, representing the upper bound for the field, is the so-called \emph{Born field strength parameter} and the nonlinear second-order operator on the left-hand side is referred to as the \emph{Born-Infeld operator}. Comparing with \eqref{MinkIntro}, we can still spot the role of the Born field strength parameter in the ratio $b/a$, while $1/a$ can be seen as a perturbation (or a diffusion) parameter, see also Section \ref{sez6}. 
\smallbreak
We incidentally recall that, in the framework of particle dynamics in the Lorentz space, an analog of \eqref{lastazionaria} (with $\beta$ replaced, roughly speaking, by $\mathfrak{c}\sqrt{m}$) originates from the relativistic Lagrangian 
$
L_r=m\mathfrak{c}^2(1-\sqrt{1-v^2/\mathfrak{c}^2}),
$
naturally arising in the construction of hypersurfaces of prescribed mean curvature in Lorentzian manifolds \cite{BarSim, CarKie}. 
Here $m$ and $v$ stand, respectively, for the mass and the velocity; the second-order operator is usually called \emph{Minkowski} (or \emph{relativistic}) operator and $\mathfrak{c}$ assumes the meaning of limit speed, equal to the speed of light. 
Of course, both the electromagnetic and the relativistic viewpoints have inspired the mathematical investigation of \eqref{Mink} and \eqref{lastazionaria}, along with the associated boundary value problems. The related bibliography is extremely rich, especially in the stationary case, see for instance \cite{Azz06, BonCoeNys, BondAvPomRei, BosColNor, CorObeOmaRiv, Maw} and the references therein. 
\smallbreak
In this article, we explain the presence of the two parameters $a, b$ in \eqref{MinkIntro} referring to the general parameter-dependent formulations of Electrostatics appearing, e.g., in \cite{Kru10, Kru2}. Thus, we interpret \eqref{MinkIntro} as a reaction-diffusion equation in 
an electromagnetic framework. On the lines of the considerations about the behavior of the Born-Infeld theory for varying $\beta$ appearing, e.g., in \cite{Gib}, the analysis of the features of the solutions of \eqref{MinkIntro} in dependence on $a$ and $b$ looks strongly motivated. The presence of the parameter $a$ allows us to carry out this investigation not only for varying field strength (Subsection \ref{sez3}), but also in a kind of vanishing/large diffusion limit (Subsection \ref{sez4}) and in dependence on a singular perturbation parameter, as well (Section \ref{sez5}).
We incidentally mention that, in literature, it is actually more common to find references to relativistic (rather than electromagnetic) reaction-diffusion models; however, we prefer to follow the described interpretation since, in principle, it may appear less justified to vary the speed of light.
\smallbreak
We will be interested in traveling waves $u(x, t)=v(x+ct)$ for \eqref{MinkIntro}, so that the \emph{wave profile} $v$ solves the second-order ODE
\begin{equation}\label{TW2}
\left(\frac{v'}{\sqrt{a^2-b^2(v')^2}}\right)' - c v' + f(v) = 0; 
\end{equation}
we will moreover require that $v(-\infty)=0$ and $v(+\infty)=1$, 
so as to deal with \emph{traveling fronts propagating with speed} $c$ for \eqref{MinkIntro}. Notice that equation \eqref{TW2} would be obtained as well in dimension $n$ when seeking \emph{planar traveling fronts}, that is, solutions of the $n$-dimensional version of \eqref{MinkIntro} having the form $u(x, t)=v(x \cdot e + ct)$, with $x \in \mathbb{R}^n$ and $e \in \mathbb{S}^1$, such that $v(-\infty)=0$ and $v(+\infty)=1$. Being \eqref{TW2} invariant with respect to translations in the independent variable, in order to recover the uniqueness of the considered fronts - so as to properly speak about their limit - we  require $v(0)=\mathcal{V}_0$ for a suitably fixed $\mathcal{V}_0 \in (0, 1)$ (see condition \eqref{normalizzazioni}). 
Since in most cases (for instance, if $f > 0$ on $(0, 1)$ and $f=0$ elsewhere) the traveling profiles are necessarily monotone (see \cite{FifMcL}), we will search from the very beginning for \emph{monotone traveling fronts}.
\smallbreak
We will see in Section \ref{sez2} (cf. also \cite{CoeSan}) that, for the families of reaction terms $f$ classically known as \emph{type A, type B and type C} \cite{BerNir}, the set of numbers $c$ for which a monotone traveling front propagating with speed $c$ exists is either made of an unbounded interval $[c^*, +\infty)$ (for type A reactions) or by a single real value $c^*$ (for type B or C reactions). In both cases the number $c^*$, depending on $a$ and $b$, takes the name of \emph{critical speed} and the corresponding traveling front solution of \eqref{MinkIntro} satisfying $v(0)=\mathcal{V}_0$ will be named \emph{critical front}. The meaning of the critical speed $c^*$ is well established in literature: in the linear diffusion case, it consists in the asymptotic speed of propagation of compactly supported initial data for the evolution equation $u_t=\Delta u + f(u)$~\cite{AroWei}.
\smallbreak
The study of front-type solutions for reaction-diffusion models with nonlinear diffusion has become quite classical in the recent literature, including a wide range of qualitatively different operators. However, the presence of models with Born-Infeld type diffusion in this context appears quite limited; as the main reference, we quote \cite{CoeSan}, where the authors investigated the traveling fronts for equation \eqref{MinkIntro} in the case $a=b=1$, mainly for $f$ positive on $(0, 1)$.
From a mathematical point of view, equation \eqref{TW2} is ruled by a particular example of \emph{singular $\phi$-Laplacian} (see, e.g., \cite{RacStaTvr06}), namely it is of the kind $(\phi(v'))'-cv'+f(v)=0$ with $\phi: ]A_-, A_+[ \to \mathbb{R}$, where $A_-, A_+ \in \mathbb{R}$, $A_- < A_+$
and $\lim_{s \to A_\pm} \vert \phi(s) \vert =+\infty$. The present investigation is carried out in the particular case $\phi(s)=s/\sqrt{a^2-b^2 s^2}$, since the possibility of performing explicit computations helps better understand the complete picture for the traveling fronts. Anyway, it is likely that a relevant part of the present results can be extended to the case of a general singular diffusion (like in \cite{CoeSan}).  As is usual when dealing with a singular operator like the Born-Infeld one, it is possible to immediately deduce some a-priori bounds on regular solutions; in particular, any $C^2$ solution of 
\eqref{Mink} necessarily satisfies $\vert \nabla u \vert < \beta$. Thus, it can be quite standard to obtain the convergence of suitable families of traveling fronts, since in case these are equibounded in $C^1$ they automatically converge locally uniformly to a Lipschitz continuous limit. 
However, sharp estimates of the critical speed for \eqref{Mink} are here less immediate to read and interpret (see \cite{CoeSan} and Proposition \ref{propoS}), due to the involved singular diffusion, and determining the shape of the limit critical profile seems far from being trivial. 
\smallbreak
In the present paper, we will proceed with the analysis of the asymptotic behavior of the critical fronts for \eqref{MinkIntro}
as $a$ and $b$ vary. More in detail, for a varying positive parameter $\gamma$, on the one hand we will deal, respectively, with the two equations
\begin{equation}\label{lista}
u_t= \left(\frac{u_x}{\sqrt{1-\gamma^2 u_x^2}}\right)_x + f(u), \qquad u_t= \gamma \left(\frac{u_x}{\sqrt{1- \gamma^2 u_x^2}}\right)_x + f(u)
\end{equation}
(Subsections \ref{sez3} and \ref{sez4}).
In the former, we focus on variations of the sole field strength ($a=1$, $b=\gamma$), while in the latter ($a=1/\gamma$, $b=1$) we replace the $1$-dimensional gradient $u_x$ by its reduced (or accelerated, depending on whether $\gamma \to 0^+$ or $\gamma \to +\infty$) counterpart $\gamma u_x$, so as to consider a kind of \emph{vanishing} (or \emph{large}) \emph{diffusion limit}. In this case, the diffusion parameter $\gamma$ and the field strength $1/\gamma$ are one the inverse of the other.  
On the other hand, in Section \ref{sez5} we will study the more delicate equation 
\begin{equation}\label{Minkintrodelta}
u_t= \gamma\left(\frac{u_x}{\sqrt{1-u_x^2}}\right)_x + f(u),
\end{equation}
corresponding to \eqref{MinkIntro} with $a=b=1/\gamma$, where we fix the field strength ($b/a=1$) and we vary the coefficient in front of the second-order operator. Such an equation may be seen, alternatively to \eqref{lista}$_{2}$, as a model for \emph{vanishing diffusion} $(\gamma \to 0^+)$ or \emph{large diffusion} ($\gamma \to +\infty)$ in the context of the considered Born-Infeld dynamics; 
however, in view of the terminology adopted for \eqref{lista}$_{2}$ and since we are mainly interested in the more difficult case $\gamma \to 0^+$, in order to avoid an unclear use of such an expression we prefer to write $\gamma=\eps$ (implying that $\eps \to 0^+$) and refer to \eqref{Minkintrodelta} as the \emph{singular perturbation problem}, leaving to the beginning of Section \ref{sez6} the considerations regarding the case $\gamma \to +\infty$. 
\smallbreak
To perform the aforementioned investigation, in Section \ref{sez2} we first provide an estimate of the critical speed associated with \eqref{MinkIntro}, complementing and partially improving the results available in literature. 
The related bounds 
will be obtained via a classical reduction of the order, for which the monotone traveling fronts for \eqref{MinkIntro} are in a one-to-one correspondence with the solutions of a suitable two-point scalar problem. 
Providing these estimates is the first step to figure out the asymptotic behaviors we are interested in; in this respect, writing $\gamma=\eps$ to indicate that $\gamma \to 0^+$ and $\gamma=1/\eps$ for $\gamma \to +\infty$, and accordingly denoting by $c_\eps^*$ the critical speed and by $v_\eps$ the associated critical profile, in Table \ref{tabella} we give a visual snapshot of the results in Sections \ref{sezz3} and \ref{sez5} for $\eps \to 0^+$. The reaction term $f$ is assumed to be continuous, with at most one sign change in $[0, 1]$, linearly controlled at $0$ and $1$, and satisfying  $\int_0^1 f(s) \, ds > 0$ (and a mild additional assumption regarding its critical points in case $a=b=1/\eps$, see Section \ref{sez5}). 
\begin{table}[ht!]
\begin{center}
\resizebox{\textwidth}{!}{\small
\begin{tabular}{|c|c|c|c|}
\hline
$b \downarrow$ $a \to$ & $\eps$  & $q$ & $1/\eps$ \\
\hline
$\eps$ &  $c_\eps^* \to +\infty$, $v_\eps \to \mathcal{V}_0$ (Sect. \ref{sez6}) & $c_\eps^* \to c^*_L$, \,\! $v_\eps \to  v_L$  (Thm. \ref{modellonuovo})
   & $c_\eps^* \, \to\,  0$, \, $v_\eps \, \to\,  H$ (Sect. \ref{sez6})   \\
\hline
$q$ & $c_\eps^* \to +\infty$, $v_\eps \to \mathcal{V}_0$  (Thm. \ref{teoremabis})  & \cellcolor{gray!70} &  $c_\eps^* \to 0$, \; $v_\eps \to H$ (Thm. \ref{conv2})   \\
\hline
$1/\eps$ & $c_\eps^* \to +\infty$,  $v_\eps \to \mathcal{V}_0$ (Sect. \ref{sez6})   & $c_\eps^* \to +\infty$, $v_\eps \to \mathcal{V}_0$ (Thm. \ref{teorema11eps}) &  
$c_\eps^* \to \bar{c} > 0$, $v_\eps \to \mathcal{V}_{L,I}$ \,\text{or}\, $\mathcal{V}_L$ (Thm. \ref{ilteorema}/Sect. \ref{sez5}) \\
\hline
\end{tabular}
}
\caption{Different behaviors of the critical speed $c_\eps^*$ and of the associated critical profile $v_\eps$ for \eqref{MinkIntro} ($\eps \to 0^+$)}
\label{tabella}
\end{center}
\end{table}
\\
Here, $q$ denotes a fixed constant (which can be taken equal to $1$, see also Remark \ref{conclusivo}); 
$c^*_L$ and $v_L$ denote, respectively, the critical speed and the associated critical front profile, satisfying $v_L(0)=\mathcal{V}_0$, for the linear diffusion equation $u_t=u_{xx} + f(u)$; $H$ is the Heaviside function, whose expression is given by $
\displaystyle H(z)=\left\{
\begin{array}{ll}
0 & z < 0 \\
1 & z > 0
\end{array}
\right.
$;
$\mathcal{V}_L: \mathbb{R} \to \mathbb{R}$ is defined as 
$
\displaystyle \mathcal{V}_L(z)=\left\{
\begin{array}{ll}
0 & z \in (-\infty, -\mathcal{V}_0)  \\
z+\mathcal{V}_0 &  z \in [-\mathcal{V}_0, 1-\mathcal{V}_0]\\
1 & z \in (1-\mathcal{V}_0, +\infty)
\end{array}
\right.
$; 
$\mathcal{V}_{L, I}$ stands for the $C^1$-gluing of $\mathcal{V}_L$ and $\mathcal{V}_I$ qualitatively depicted in Figure \ref{introduttiva} for illustrative purposes (the formal details being postponed to Section \ref{sez5}), where
$\mathcal{V}_I: \mathbb{R} \to \mathbb{R}$ is the solution of the inviscid problem 
$\displaystyle 
\left\{
\begin{array}{l}
\bar{c} v' = f(v) \\
v(0)=\mathcal{V}_0
\end{array}
\right.
$,
being $\bar{c}=\lim_{\eps \to 0^+} c_\eps^*$.
We notice that the indicated 
convergences hold \emph{independently of the considered reaction terms}; in this respect, we underline that the precise information about the asymptotic behavior of $c_\eps^*$ for $\eps \to 0^+$ in case $f$ is not strictly positive in a neighborhood of $0$ (that is, for type B or type C reactions, see Section \ref{sez2}) can be detected in a precise way thanks to the properties of the Born-Infeld operator, while it may not be the case for linear or saturating (cf. \cite{Ros}) diffusions. 
With these preliminaries, we briefly explore the results in Table \ref{tabella}:
\begin{itemize}
\item first, quite naturally, if we strengthen the singularity sending the maximal field strength to $0$ - namely, $b \to +\infty$ - while $a$ is kept fixed (or $a \to 0^+$, keeping $b$ fixed), then the critical profiles have no room for increasing their derivative, so they necessarily converge to the constant value taken at $z=0$; on the contrary, if the singularity is weakened by sending $a \to +\infty$ while $b$ is kept fixed, then even jumping limit profiles are admitted, in spite of the singular diffusion (since the bounds on the derivative of $v_\eps$ become arbitrarily large);
\item a little subtler is the case in which $a > 0$ is constant and $b \to 0^+$ - apparently analogous to $a \to +\infty$, $b$ fixed - where in the limit we do not find a stepwise profile, but rather we are able to rigorously recover the critical profile of the linear diffusion case. This perfectly agrees with the fact that, quoting \cite[p. 4]{CarKie}, 
\begin{center}
\!\!\!\!\textbf{in the limit $\beta \to 0$, Born's law goes over into Maxwell's law of the ``pure aether''}.
\end{center}
In other words, the picture for $b \to 0^+$ is in line with the physical predictions according to which the classical theory is recovered by admitting arbitrarily large electric fields \cite{Gib};
\item finally, for the singular perturbation problem  ($a=1/\eps, b=1/\eps$) with $f$ of type A or B, we obtain a completely new and rather surprising result, which is actually the main finding of the present work: for $\eps \to 0^+$, 
\begin{center}
\textbf{the critical speed $c_\eps^*$ converges to a positive value $\bar{c}$ and, in the limit,}  
\textbf{the critical front profile becomes sharp on one side only, near the value~$0$ ($v_\eps \to \mathcal{V}_{L, I}$)}. 
\end{center}
\end{itemize}
As for this last result, in Section \ref{sez5} we will explicitly show that $\bar{c}=f(v_+)$, 
where $v_+$ is the largest root of the equation $F(v)=vf(v)$ (being $F(v)=\int_0^v f(s) \, ds$). 
It has to be said that the natural limit profile one could expect for the Born-Infeld equation in presence of a small perturbation parameter would be the fully piecewise linear one $\mathcal{V}_L$ (see, e.g., \cite{BosFel}). In fact, it is easy to see (Subsection \ref{ibilanciati}) that this is the case for balanced reaction terms of type C, for which $c_\eps^* = 0$ for every $\eps$ (more in general, in Section \ref{sez5} we will see that if $v_\eps \to \mathcal{V}_L$, then necessarily $\bar{c}=\int_0^1 f(s) \, ds$). To the author's knowledge, the fact that $\bar{c} > 0$ and the onset of a solely one-sided sharpening of the limit front profile are instead new phenomena in the theory of traveling waves with nondegenerate diffusion, not included in the classical literature regarding sharp fronts (mainly for the degenerate case) as, e.g., \cite{MalMar03, SanMai94}. We notice as well the profound differences with the usual stepwise limit appearing for linear \cite{HilKim} and saturating \cite{Gar} diffusions, for $\bar{c}=0$.
\smallbreak
\noindent
The conclusions of Sections \ref{sezz3}-\ref{sez5} are complemented by the discussion in Section \ref{sez6}, where we also provide a further snapshot of the limit picture for equation \eqref{TW2} for varying $a, b$ (Figure \ref{comportamenti}), and are supported by several numerical experiments performed with Wolfram Mathematica$^\copyright$ software (see also the brief Appendix concluding the paper).

\section{An elementary estimate of the critical speed}\label{sez2}

In this section, we provide an explicit estimate of the critical speed associated with \eqref{TW2}, based on a simple observation on the associated first-order reduction. 
\smallbreak
To this end, henceforth we assume that 
\begin{center}
(H) \; $f \in C([0, 1])$ satisfies $f(0)=f(1)=0$ and 
there exists $k > 0$ for which \\ $\vert f(s) \vert \leq ks$, $\vert f(s) \vert \leq k(1-s)$ for every $s \in [0, 1]$.
\end{center} 
The linear controls at $0$ and at $1$ are very common in the theory of traveling waves and ensure that the found wave fronts are proper, that is, defined on the whole $\mathbb{R}$ (see, e.g., \cite{GarSan} and \cite[Theorem 10.5]{GilKer}); however, we are not requiring that $f'(0)$ exists, as is quite usual in this kind of models. 
Among the reaction terms satisfying (H), the most common ones are the following:
\begin{itemize}
\item[(A)] $f$ of type A (or \emph{monostable}), i.e., $f(u) > 0$ for every $u \in (0, 1)$;
\item[(B)] $f$ of type B (or \emph{of combustion type}), i.e., there exists $\alpha \in (0, 1)$ for which $f(u)=0$ for $u \in [0, \alpha]$ and $f(u) > 0$ for every $u \in (\alpha, 1)$;
\item[(C)] $f$ of type C (or \emph{bistable}), i.e., there exists $\alpha \in (0, 1)$ for which $f(s)<0$ for $s \in (0, \alpha)$, $f(s) > 0$ for $s \in (\alpha, 1)$. In this case, we assume $\int_0^1 f(s) \, ds \geq 0$. 
\end{itemize}
We are concerned with monotone wave fronts $u(x, t)=v(x+ct)$ for \eqref{MinkIntro}, corresponding to \emph{monotone heteroclinic solutions of \eqref{TW2} joining $0$ and $1$} - henceforth, briefly called \emph{heteroclinics}. The usual reduction of the order going back to \cite{FifMcL} and widely exploited in literature (see, e.g., \cite{GarSan, MalMar}), based on the possibility of inverting the monotone function $v=v(z)$, obtaining then a first-order differential equation satisfied by $\phi(v)=v'(z(v))$, leads naturally to the problem
\begin{equation}\label{notazioney}
\left\{
\begin{array}{l}
\displaystyle y'= c a \frac{\sqrt{y(2a+b^2y)}}{a+b^2y} - f(v) \vspace{0.1cm}\\
y(0)=0, \, y(1)= 0, \, y > 0 \text{ on } (0, 1),
\end{array}
\right.
\end{equation}
where 
\begin{equation}\label{costruzionev}
y(v)=\frac{1}{b^2} \left(\frac{a^2}{\sqrt{a^2-b^2 v'(z(v))^2}}-a\right) \quad \text{ and } \quad 
v'(z)=\frac{a\sqrt{y(v(z))(2a+b^2 y(v(z))}}{a+b^2 y(v(z))};
\end{equation}
notice that the requirements $v(-\infty)=0$ and $v(+\infty)=1$ give rise to the two boundary conditions in \eqref{notazioney}.
\begin{definition}
We say that $c \in \mathbb{R}$ is an \emph{admissible speed} if problem \eqref{notazioney} has a solution. 
\end{definition}
In case $c \in \mathbb{R}$ is admissible, thanks to assumption (H) it can be seen that the differential equation in \eqref{costruzionev} allows one to reconstruct a strictly monotone front profile $v$ defined on the whole $\mathbb{R}$ and joining, at $\pm \infty$, the equilibria $0$ and $1$. Multiplying \eqref{TW2} by $v'$ and integrating on the whole real line, moreover, it can be seen that $c$ has the same sign as $\int_0^1 f(s) \, ds$, so throughout the present work we will always deal with \emph{positive} admissible speeds, in view of the above assumptions. 
\\
For further use, we set 
$$
R_{a, b}(s):=\frac{\sqrt{s(2a+b^2s)}}{a+b^2s}
$$
and we notice that $R_{a, b}$ is \emph{increasing and bounded}, fulfilling the estimate $R_{a, b}(s) \leq 1/b$ for every $s > 0$. The boundedness of $R_{a, b}$ is a consequence of the presence of a singular $\phi$-Laplacian in our model (see the Introduction) and is the main difference with respect to the linear (and also with the saturating, see \cite{Gar}) diffusion case; it will turn into qualitatively different results about the asymptotic behavior of front profiles, as already mentioned in the Introduction. 
\\
Problem \eqref{notazioney} is usually dealt with by studying separately 
the two Cauchy problems (respectively, forward and backward)
\begin{equation}\label{cauchyf}
(P_{c, f})^+ \quad \left\{
\begin{array}{l}
\displaystyle y'= ca R_{a, b}(y) - f(v) \vspace{0.1cm}\\
y(0)=0 \,
\end{array}
\right.
\quad \text{and} \quad  (P_{c, f})^- \quad \left\{
\begin{array}{l}
\displaystyle y'= ca R_{a, b}(y) - f(v) \vspace{0.1cm}\\
y(1)=0, \,
\end{array}
\right.
\end{equation}
whose solutions are denoted by $y_{c, f}^+$ and $y_{c, f}^-$, respectively. 
If $c$ is an admissible speed, then $y_{c, f}^+ \equiv y_{c, f}^-$; to this end, by the uniqueness, it is sufficient that there exists $v_0 \in (0, 1)$ for which $y_{c, f}^+(v_0)=y_{c, f}^-(v_0)> 0$. One then usually attempts to deduce some bounds on the admissible speeds via the analysis of the mutual position of the graphs of $y_{c, f}^+$ and $y_{c, f}^-$ as $c$ varies, often invoking continuous dependence and monotonicity arguments for $(P_{c, f})^-$, which enjoys uniqueness due to the monotonicity of $R_{a, b}$. \\
So far, to the author's knowledge, the best estimate of the admissible speeds for \eqref{notazioney} obtained in this way has been provided in \cite{CoeSan}, in the case $a=b=1$, as follows; notice that the argument therein works as well for a general second-order equation ruled by a singular $\phi$-Laplacian.
\begin{theorem}\cite[Proposition 3.2]{CoeSan}\label{propoS}
Let $a=b=1$ and assume that $f$ is a type A function satisfying (H). Moreover, for some number $M > 0$, let 
\begin{equation}\label{controllomink}
f(s) \leq \frac{Ms}{\sqrt{1+M s^2}} \quad \textrm{ for every } s \in [0, 1].
\end{equation}
Then, the set of the admissible speeds for \eqref{TW2} is an interval $[c^*, +\infty)$, where $c^* \leq 2\sqrt{M}$. 
\end{theorem}

Due to the way the result is obtained, it is here convenient to control the growth of $f$ through the function $\mathcal{B}_M(s):=Ms/\sqrt{1+Ms^2}$, a little less immediate to interpret than the usual linear control appearing in (H). Moreover, \eqref{controllomink} does not rescale homogeneously upon dilation of $f$ (in the sense that if $f$ satisfies \eqref{controllomink}, for small $\eps$ the function $g(s)=f(s)/\eps$ does not fulfill \eqref{controllomink} with $M$ replaced by $M/\eps$), and this could make the asymptotics for the critical speed less clear. 
\smallbreak
\noindent
In the following proposition, we provide new estimates of the critical speed which complement and improve Theorem \ref{propoS} in some cases (see Remark \ref{rconfronto}), for general reaction terms belonging to the aforementioned classes. Henceforth, we set $F(v)=\int_0^v f(s) \, ds$.  
\begin{proposition}\label{proponuova}
Assume that $f \in C([0, 1])$ satisfies (H).
Then, the following hold: 
\begin{enumerate}
\item if $f$ is of type A, then the set of the admissible speeds for problem \eqref{notazioney} is an interval $[c^*, +\infty)$, where 
\begin{equation}\label{Stimaa}
\sqrt{\sup_{v \in (0, 1]}\left(\frac{b^2}{a^2}\frac{F(v)^2}{v^2}+\frac{2}{a}\frac{F(v)}{v^2}\right)}
\leq c^* \leq \sup_{v \in (0, 1]} \left(\frac{b}{a} f(v) + 2\sqrt{\frac{f(v)}{av}} \right). 
\end{equation}
\item if $f$ is of type B or C, then 
there exists a unique admissible speed $c^*$ satisfying the estimate 
\begin{equation}\label{Stimac}
\frac{1}{R_{a, b}(F(1)-F(\alpha))}  \sup_{v \in (0, 1]}\frac{F(v)}{av} \leq  c^*  \leq \frac{1}{\alpha} \sqrt{\frac{b^2}{a^2} (F(1)-F(\alpha))^2+\frac{2}{a}(F(1)-F(\alpha))}.
\end{equation}
\end{enumerate}
\end{proposition}

\begin{proof}
\begin{enumerate}
\item The fact that the admissible speeds form a nonempty unbounded interval $[c^*, +\infty)$ can be proved as in \cite[Proposition 3.2]{CoeSan}. We thus prove the bounds on $c^*$ appearing in \eqref{Stimaa}. 
As for the lower bound, we notice that if $c$ is an admissible speed, then the corresponding solution $y$ of \eqref{notazioney} satisfies 
$$
\left\{
\begin{array}{l}
\displaystyle y'= c a R_{a, b}(y) - f(v) \leq c a \frac{\sqrt{y(2a+b^2 y)}}{a+b^2 y} \vspace{0.1cm}\\
y(0)=0. 
\end{array}
\right.
$$
Denote by $y_m$ the maximal solution (notice that here there is no uniqueness) of $(P_{c, 0})^+$, explicitly given by 
$y_m(v)=\frac{a}{b^2} (\sqrt{1+c^2 b^2 v^2}-1)$.
Standard lower and upper solution theory for first-order ODEs (see, e.g., \cite{Wal}) then provides $y(v) \leq y_m(v)$ for every $v \in [0, 1]$. 
From the monotonicity of $R_{a, b}$, we can thus infer that 
$$
y'(v) \leq ca R_{a, b}(y_m(v)) - f(v)  = \frac{c^2 a v}{\sqrt{1+c^2b^2v^2}} - f(v) \; \Rightarrow \; y(v) \leq \frac{a}{b^2}(\sqrt{1+c^2 b^2 v^2}-1)-F(v), 
$$
where the latter inequality is obtained integrating the former one between $0$ and $v$. Since the last expression has to be positive for every $v \in (0, 1)$, the conclusion follows. 
\\
We now prove the upper bound in \eqref{Stimaa}; to this end, similarly as in \cite{CoeSan}, we seek a positive lower solution of $(P_{c, f})^+$ by solving problem $(P_{\beta, 0}^+)$ in dependence on a positive parameter $\beta$, to be found in such a way that the maximal (positive) solution 
$
y_{\beta, 0}^+(v)=\frac{a}{b^2} (\sqrt{1+\beta^2 b^2 v^2}-1)
$
of $(P_{\beta, 0}^+)$
solves the differential inequality
$$
y'(v) \leq ca R_{a, b}(y(v)) - f(v).
$$
Replacing the explicit expression of $y_{\beta, 0}^+$, we thus have to find $\beta$ for which 
\begin{equation}\label{disugbeta}
\frac{\beta^2 a v}{\sqrt{1+\beta^2 b^2 v^2}}-\frac{c\beta a v}{\sqrt{1+\beta^2 b^2 v^2}}+ f(v) \leq 0 
\end{equation}
for every $v \in (0, 1]$ (notice that \eqref{disugbeta} is satisfied for $v=0$ in view of (H)); 
dividing all the summands by $av$ and using the fact that $\sqrt{1+b^2 \beta^2 v^2} \leq 1+b\beta v$, this will be true if for every $v \in (0, 1]$ it holds
\begin{equation}\label{relazione}
\beta^2-\left(c-\frac{b f(v)}{a}\right)\beta+ \frac{f(v)}{av}\leq 0;
\end{equation}
in order for such a $\beta$ to exist, $c$ has then to satisfy
\begin{equation}\label{ultima}
\left(c-\frac{b f(v)}{a}\right)^2-4\frac{f(v)}{av} \geq 0 \quad \Rightarrow \quad c \geq \frac{b}{a} f(v) + 2 \sqrt{\frac{f(v)}{av}}  
\end{equation}
for every $v \in (0, 1]$.
The thesis follows.
\item 
After noticing that, in this second case, uniqueness holds also for $(P_{c, f})^+$ (see for instance \cite{BonSan}), the existence and the uniqueness of an admissible speed follow from a standard continuity and monotonicity argument. Precisely, on the one hand $y_{0, f}^+(\alpha) < y_{0, f}^-(\alpha)$, since $F(1)=\int_0^1 f(s) \, ds > 0$. On the other hand, taking into account that $y_{0, f}^-(v)=\int_v^1 f(s) \, ds$ is an upper solution for $(P_{c, f})^-$ (recall that $c$ is positive), we have that $y_{c, f}^-(\alpha) \leq \int_\alpha^1 f(s) \, ds=F(1)-F(\alpha)$ for every $c > 0$, whereas $y_{c, f}^+(\alpha) > y_{c, 0}^+(\alpha) \to +\infty$ for $c \to +\infty$. Being $y_{c, f}^+(\alpha), y_{c, f}^-(\alpha)$ continuous and strictly monotone with respect to $c$, necessarily there exists a (unique, by strict monotonicity) value $c^*$ for which $y_{c^*, f}^+(\alpha)=y_{c^*, f}^-(\alpha)$. The claim follows. 
\\
As for the lower bound on $c^*$, 
the fact that $y_{0, f}^-(v)$ is an upper solution for $(P_{c, f})^-$ on the whole interval $[0, 1]$, for every positive $c$, implies that $\Vert y_{c^*, f}^+ \Vert_{L^\infty(0, 1)} \leq \Vert y_{0, f}^- \Vert_{L^\infty(0, 1)}=\int_\alpha^1 f(s) \, ds=F(1)-F(\alpha)$. 
Since $R_{a, b}$ is increasing, 
the solution of \eqref{notazioney} then satisfies
$$
y'(v) \leq c^* a R_{a, b}(F(1)-F(\alpha))  -f(v)
$$
for every $v \in [0, 1]$. Integrating this last relation from $0$ to $v$ provides the desired bound, taking into account that $y(v)$ has to be positive for every $v \in (0, 1)$ and $y(0)=0$.
\\ 
Finally, to prove the upper bound on $c^*$ we observe that, being 
$$
y_{c^*, f}^+(\alpha) \geq \frac{a}{b^2}(\sqrt{1+(c^*)^2b^2 \alpha^2}-1), \qquad y_{c^*, f}^-(\alpha) \leq y_{0, f}^-(\alpha) = F(1)-F(\alpha),
$$
it has necessarily to hold $\dfrac{a}{b^2}(\sqrt{1+(c^*)^2b^2 \alpha^2}-1) \leq F(1)-F(\alpha)$ in order for 
$y_{c^*, f}^+(\alpha)=y_{c^*, f}^-(\alpha)$ to be fulfilled. The thesis follows.
\end{enumerate}
\end{proof}
Observe
that the procedure employed to find the lower bounds in \eqref{Stimaa} and \eqref{Stimac}, based on determining a global estimate for $y$ to be inserted into the differential equation in \eqref{notazioney}, is used for instance in \cite{DraTak} in the context of traveling waves. A simple way of providing such a global estimate is through the upper solution $y_{0, f}^-$, as in the proof of Item 2 of Proposition \ref{proponuova} (this would work also in the type A case, yielding however a less sharp lower bound for $c^*$ than the one in \eqref{Stimaa}); this is independent of the considered diffusion, provided that the term replacing $R_{a, b}$ in the differential equation for $y$ is increasing. For instance, the same procedure may be exploited in the linear diffusion case, where the differential equation governing the first-order reduction reads as $y'=2c\sqrt{y}-2f(v)$.
One could also think about bootstrapping such a scheme, finding a sequence of sharper (but cumbersome) bounds for $y$ and hence for $c^*$; however, it seems that this does not improve significantly the asymptotic trend of the critical speed. 
Incidentally, note as well that the proof of \eqref{Stimac} could be slightly refined, invoking the validity of the inequalities not only at $v=\alpha$, but along a whole interval. This would lead to finer bounds, however it would again penalize the readability of \eqref{Stimac} without providing any substantial improvement, reason for which we have not pursued this possibility. 
Finally, we notice that \eqref{Stimaa} and \eqref{Stimac} imply the simpler estimate
\begin{equation}\label{forseutile}
c^* \geq \dfrac{b}{a} \sup_{v \in (0, 1]} \dfrac{F(v)}{v},
\end{equation} 
which would directly follow using the inequality $R_{a, b} \leq 1/b$ in the differential equation for $y$ and integrating, similarly as before. It is worth mentioning that it is the presence of the singular $\phi$-Laplacian, via the consequent boundedness of $R_{a, b}$, that allows one to infer such a ``universal'' lower bound, which is \emph{independent of the shape of $f$}. 
The same holds indeed for any first-order reduction of the kind $y'=cH(y) - f(v)$, with $H$ a bounded function.

\begin{definition}
The quantity $c^*$ appearing in the statement of Proposition \ref{proponuova} is called \emph{critical speed} and the corresponding heteroclinic solution of \eqref{TW2} increasingly connecting $0$ and $1$, unique up to $z$-translations, is referred to as the \emph{critical front profile} (or simply \emph{critical profile}).
\end{definition}
\begin{corollary}
Let $a=b=1$. Then, keeping the same notation as in Proposition \ref{proponuova}, the following hold: 
\begin{enumerate}
\item if $f$ is of type A, then the critical speed $c^*$ satisfies the estimate
$$
\sqrt{\sup_{v \in (0, 1]}\left(\frac{F(v)^2}{v^2}+\frac{2F(v)}{v^2}\right)} \leq c^* \leq \sup_{v \in (0, 1]} \left(f(v)+2\sqrt{\frac{f(v)}{v}} \right);
$$
\item if $f$ is of type B or C, then the critical speed $c^*$ satisfies the estimate
$$
\frac{1+F(1)-F(\alpha)}{\sqrt{(F(1)-F(\alpha))(2+F(1)-F(\alpha))}} \sup_{v \in (0, 1]}\!\! \frac{F(v)}{v} \! \leq \!  c^*  \!\leq \! \frac{1}{\alpha}\sqrt{(F(1)-F(\alpha))^2\!+\!2(F(1)-F(\alpha))}.
$$
\end{enumerate}
\end{corollary}
We close this section stating a result concerning the monotonicity of the critical speed with respect to the parameters $a$ and $b$, which will turn to be useful henceforth (see also \cite[Proposition 2.8]{Gar}). To this end, we denote the critical speed by $c_{a, b}^*$, to better highlight its dependence on the parameters.
\begin{proposition}\label{monotono}
Let $f$ be of type A, B or C (satisfying, in particular, assumption (H)). Then:
\begin{itemize}
\item[-] for fixed $b > 0$, $a \mapsto c_{a, b}^*$ is monotone decreasing;
\item[-] for fixed $a > 0$, $b \mapsto c_{a, b}^*$ is monotone increasing. 
\end{itemize}
\end{proposition}
\begin{proof}
We prove the first statement. Let $\bar{a} > 0$ be fixed and $c_{\bar{a}, b}^*$ be the corresponding critical speed. 
Since, for every $y \geq 0$, $a \mapsto R_{a, b}(y)$ is monotone increasing, if $a > \bar{a}$ then for every $v \in [0, 1]$ one has $y_{c_{a, b}^*, f}^- (v) \leq y_{c_{\bar{a}, b}^*, f}^- (v)$ for every $v \in [0, 1]$. It follows that either $c_{\bar{a}, b}^*$ is admissible for \eqref{notazioney} (and hence $c_{a, b}^* \leq c_{\bar{a}, b}^*$) or $y_{c_{a, b}^*, f}^-(v_0) = 0$ for some $v_0 > 0$; in this second case, one has to take $c < c_{a, b}^*$ in order for $y_{c, f}^-$ to fulfill \eqref{notazioney}, hence $c_{a, b}^* < c_{\bar{a}, b}^*$.
As for the second statement, the argument may be repeated taking into account that $b \mapsto R_{a, b}(y)$ is monotone decreasing.
\end{proof}
The behavior highlighted by Proposition \ref{monotono} is also suggested by the dependences on $a$ and $b$ appearing in \eqref{Stimaa} and \eqref{Stimac}.

\subsection{Some further remarks}

\begin{remark}\label{KPPrem}
\textnormal{\emph{(Behavior of $f$ at $0$ and admissible speeds: complementing and further discussing \eqref{Stimaa}).}
In the particular case when $f$ is of type A and $f'(0)$ exists, one can further deduce that
$c^* \geq 2\sqrt{f'(0)/a}$ by proceeding essentially as in \cite{CoeSan} and \cite[Lemma 3.1]{GarSan}. For a type A function $f$ satisfying 
\begin{equation}\label{KPP}
f(s) \leq f'(0)s \quad \text{ for every } s \in [0, 1],
\end{equation}
namely one of the so-called \emph{KPP reaction terms}, one can then infer an estimate of the kind
$$
\max\left\{2\sqrt{\frac{f'(0)}{a}}, \sqrt{\frac{b^2F(1)^2}{a^2}+\frac{2F(1)}{a}} \right\} \leq c^* \leq 2\sqrt{\frac{f'(0)}{a}}+\frac{b\Vert f \Vert_\infty}{a};
$$
of course, this says, for instance, that the critical speed converges to $2 \sqrt{f'(0)/a}$ in case $b \to 0^+$ (see Section \ref{sezz3}). However, we underline that the estimates obtained in Proposition \ref{proponuova} hold and provide an explicit strictly positive lower bound on $c^*$ also in the ``dangerous'' cases when $f'(0)$ exists but $f'(0)=0$ - in which the direct analysis of the equilibrium $(0, 0)$ in the phase plane would only ensure $c^* \geq 0$ - or when $f'(0)$ does not exist, as already remarked. In this latter case, one could infer further information if both $\ell:=\liminf_{s \to 0^+} f(s)/s$ and $\limsup_{s \to 0^+} f(s)/s$ belong to $\mathbb{R}$, by observing (analogously as in \cite[Proposition 3.4]{CoeSan}) that $E(y)=\int_0^y ds/R_{a, b}(s) =\sqrt{y(2a+b^2y)}$ satisfies
$$
E(y(v))'=ca -\frac{f(v)}{R_{a, b}(y(v))}, 
$$
whence $L:=\limsup_{v \to 0^+} E(y(v)) \geq 0$ is such that $L^2-caL+\ell a \leq 0$. Consequently, one could deduce $c^2 a^2 -4 \ell a \geq 0$, which provides a significant lower bound on $c^*$ if $\ell >0$ (thus necessarily in the type A case). \\
We finally observe that if $f$ satisfies instead $\lim_{s \to 0^+} f(s)/s = +\infty$ (so that it has infinite slope at $0$ and (H) does not hold) and is positive on $(0, 1)$, by proceeding as in \cite[Proposition 4.2]{DraTak} it can be shown that \emph{no admissible speeds exist}. Indeed, using $y(v) \leq y_m(v)=\frac{a}{b^2} (\sqrt{1+c^2 b^2 v^2}-1)$  
in the differential equation for $y$ would yield $(\sqrt{y(v)(2a+b^2y(v))})' \leq ca - \frac{f(v)}{cav} (a+b^2 y(v)) \leq ca - \frac{f(v)}{cv}$, leading to the contradiction $y \equiv 0$.  
} 
\end{remark}

\begin{remark}\label{largeroot}
\textnormal{\emph{(Characterization of the critical fronts).} In case $f'(0)$ exists, we have a further piece of information about the critical fronts: thanks to \cite[Propositions 3.4 and 4.2]{CoeSan}, adapted to our setting which includes the parameters $a$ and $b$, we infer that the funcion $E(y)$ introduced in Remark \ref{KPPrem} is such that 
$E(y(v))'(0)$ is the largest root of the equation $x^2-cax+af'(0)=0$. Notice that this is completely \emph{independent of} $b$. The critical speed is thus characterized by the asymptotic behavior of the traveling profile at $z=-\infty$, as in the linear case. 
}
\end{remark}

\begin{remark}\label{rconfronto}
\textnormal{\emph{(The upper bound on $c^*$: comparison with the estimates in \cite{CoeSan}).}
Assume $f$ of type A and $a=b=1$, and recall that we have set $\mathcal{B}_M(s)=Ms/\sqrt{1+Ms^2}$. We first notice that if $v_0 \in (0, 1)$ is such that $f(v_0)=f_M:=\max_{v \in [0, 1]}f(v)$, from $f_M \leq \mathcal{B}_M(v_0)$ and the fact that $\mathcal{B}_M$ is increasing it turns out that 
the least value of $M$ for which \eqref{controllomink} may hold produces
\begin{equation}\label{stima1}
c^* \leq \displaystyle \sqrt{\frac{2f_M(f_Mv_0+\sqrt{4+f_M^2 v_0^2})}{v_0}}.
\end{equation}
Whether \eqref{stima1} is sharper than \eqref{Stimaa} depends of course on $f_M$ and $v_0$; just to give an idea, in case of a Fisher-type reaction term $f(s)=ms(1-s)$, $m > 0$, the leading terms in \eqref{stima1} for $m \to 0$ and $m \to +\infty$ are $\sqrt{2m}$ and $m/2$, respectively, against the counterparts $2\sqrt{m}$ and $m/4$ obtained considering separately each summand of the upper bound in \eqref{Stimaa} (this choice is not optimal, but allows one to maintain a certain readability). Hence, \eqref{Stimaa} is sharper for $m$ sufficiently large (it can be seen that it starts becoming sharper for $m \approx 31.9$). 
For instance, it may be easily checked that \eqref{Stimaa} improves Theorem \ref{propoS} for
$$
f_1(s)=6s(1-s)(1+6s), \quad f_2(s)=40s^2(1-s); 
$$
the former is the classical positive nonlinearity $f(s)=ms(1-s)(1+\sigma s)$ dealt with in \cite{AroWei, HadRot}, while the latter is a Huxley-type \cite{HodHux} nonlinearity (but the same could be done with a more general Nagylaki-type \cite{Nag75} nonlinearity $f(s)=ms(1-s)(1+\sigma-2\sigma s)$, for $\sigma \approx -1$). 
More in general, the higher $v_0$, the more the behavior of the function $\mathcal{B}_M$ produces a worse estimate for large $f_M$, since the optimal choice of $M$ approaches the value $M/\sqrt{1+M}$ and it has to be $M/\sqrt{1+M} \geq f_M$. The upper bound in \eqref{Stimaa} has instead the advantage of being more immediately readable, also in presence of rescaling parameters, and more ``symmetric'', not being affected by the region where $f$ reaches its maximum. Proposition \ref{proponuova} thus
complements the results in \cite{CoeSan}, possibly improving them for large reaction terms, and is suitable for our considerations on the limit critical profiles in Sections \ref{sezz3} and \ref{sez5}. Nevertheless, we could also proceed as in \cite[Proposition 3.2]{CoeSan} in order to obtain the following similar result. } 
\begin{proposition}
Let $f$ satisfy (H) and (A). Moreover, for some number $M > 0$, let 
\begin{equation}\label{controllomink2}
f(s) \leq \frac{Ms}{\sqrt{1+\frac{M}{a}b^2 s^2}} \quad \textrm{ for every } s \in [0, 1].
\end{equation}
Then, the set of the admissible speeds for \eqref{TW2} is an interval $[c^*, +\infty)$, where $c^* \leq 2\sqrt{\frac{M}{a}}$. 
\end{proposition}
\textnormal{
The statement can be obtained by mimicking the argument in the proof of \cite[Proposition 3.2]{CoeSan}, seeking $\beta > 0$ such that $y_{\beta, 0}^+$ satisfies the inequality
$
\displaystyle y'(v) \leq ca R_{a, b}(y(v)) - \frac{Mv}{\sqrt{1+\beta^2 b^2 v^2}},
$
in view of \eqref{controllomink2}. One then exploits the fact that $\beta$ can be chosen less than or equal to $\sqrt{\frac{M}{a}}$ to conclude. In general, however, this estimate does not improve in a relevant way the asymptotic behavior of $c^*$ with respect to the parameters $a, b$. } 
\end{remark} 

\begin{remark}\label{finoaK}
\textnormal{\emph{(Other improvements of the upper bound in \eqref{Stimaa}).}} 
\textnormal{
In the proof of Item 1 of Proposition \ref{proponuova}, in principle it is not necessary that $y_{\beta, 0}^+$ is a lower solution for $(P_{c, f})^+$ on the whole interval $[0, 1]$, but it could suffice that this occurs  
only on a certain interval $[0, K]$, $K < 1$. To this end, similarly as in \eqref{ultima}, one would require
\begin{equation}\label{velocparziale}
c \geq \sup_{v \in (0, K]} \left(\frac{b}{a} f(v)+ 2 \sqrt{\frac{f(v)}{av}}\right).
\end{equation}
Taking $c=\sup_{v \in (0, K]} \left(\frac{b}{a} f(v)+ 2 \sqrt{\frac{f(v)}{av}}\right)$,
if one can find a suitable $\beta > 0$ for which \eqref{relazione} holds for every $v \in [0, K]$, it will be
$y_{c, f}^+(K) \geq y_{\beta, 0}^+(K)=\frac{a}{b^2}(\sqrt{1+\beta^2 b^2 K^2}-1)$. Now, the solution of 
$$
\left\{
\begin{array}{l}
\displaystyle y'= -f(v) \vspace{0.1cm}\\
y(K)=\displaystyle \frac{a}{b^2}(\sqrt{1+\beta^2 b^2 K^2}-1),
\end{array}
\right.
$$
given by $\hat{y}(v)=\displaystyle \frac{a}{b^2}(\sqrt{1+\beta^2 b^2 K^2}-1)-\int_K^v f(s) \, ds$, is a lower solution for $(P_{c, f})^+$ in the interval $[K, 1]$, so that, if $\hat{y}(1) \geq 0$, then
$$
y_{c, f}^+(v) \geq \tilde{y}(v):=
\left\{
\begin{array}{ll}
y_{\beta, 0}^+(v) & v \in [0, K] \\
\hat{y}(v) & v \in (K, 1]
\end{array}
\right. \quad \text{ for every } v \in [0, 1].
$$ 
For instance, if the right-hand side in \eqref{velocparziale} equals $\sup_{v \in (0, K]} \frac{b}{a} f(v) + \sup_{v \in (0, K]}2 \sqrt{\frac{f(v)}{av}}$, one can choose $\beta=\sup_{v \in (0, K]}\sqrt{\frac{f(v)}{av}}$ and, if $\hat{y}(1) \geq 0$, it follows that 
$$
c^* \leq \sup_{v \in (0, K]} \left(\frac{b}{a} f(v) + 2\sqrt{\frac{f(v)}{av}}\right).
$$
Such an argument is potentially applicable also to other kinds of operators. As another refinement of \eqref{Stimaa}, moreover, we observe that the explicit expression of the best constant $\beta$ for which \eqref{disugbeta} is fulfilled would actually be obtainable (even if slightly involved), possibly producing a finer upper bound for the critical speed. However, none of these improvements would provide further relevant information about the asymptotic behavior of $c^*$ in dependence on $a$ and $b$.
}
\end{remark}

\begin{remark}\label{bilanciata}
\textnormal{\emph{($0$-critical speed).}
Under our assumptions, the only possibility to have zero critical speed, yielding a \emph{stationary wave front}, is for a reaction term of type C satisfying $\int_0^1 f(s) \, ds = 0$ (usually referred to as \emph{balanced} reaction). In this case, with reference to \eqref{notazioney}, it is immediate to check that $c=0$ is indeed the critical speed. }
\end{remark}

\begin{remark}
\textnormal{\emph{(Fronts connecting $\alpha$ and $1$).}
Clearly, in the type B and type C cases, one could give estimates of the critical speed analogous to \eqref{Stimaa} for traveling fronts connecting the equilibria $\alpha$ and $1$, provided that assumption (H) holds also at the equilibrium $\alpha$; indeed, $f$ restricted to $[\alpha, 1]$ becomes a type A function. In the same spirit, one could also describe the asymptotic behavior of the corresponding critical profiles, as in the next sections (and similarly as in \cite{Gar}, for instance). We omit the explicit statements for briefness. 
}
\end{remark}

\section{Some physically predicted convergence results}\label{sezz3}

We here give some results which are in line with the physical predictions about the Born-Infeld model, first focusing on variations of the field strength (Subsection \ref{sez3}) and then considering the vanishing/large diffusion limit (Subsection \ref{sez4}).

\subsection{Behavior of the critical profile for varying field strength}\label{sez3}

In this section, we analyze the behavior of the critical front profile for \eqref{MinkIntro} according to variations of the field strength parameter $1/b$, assuming $a$ to be constant. Without loss of generality, we can take $a=1$ so that, setting $b=\gamma$, we consider 
the equation 
\begin{equation}\label{eqBI}
\left(\frac{v'}{\sqrt{1-\gamma^2 (v')^2}}\right)' - c_\gamma^* v' + f(v) =0.
\end{equation}
Here we have already made the \emph{critical speed} $c_\gamma^*$ explicit, by virtue of the results in Section \ref{sez2}.
In order to recover the uniqueness of the critical traveling profile $v_\gamma$, henceforth we will impose that it satisfies the ``normalization'' condition
\begin{equation}\label{normalizzazioni}
v_\gamma(0)=\mathcal{V}_0, \; \text{ where } \;\mathcal{V}_0=\left\{\begin{array}{ll} 1/2 & \text{ if } f \text{ is of type A} \\ \alpha & \text{ if } f \text { is of type B or C,} \end{array} \right.
\end{equation}
unless explicitly stated otherwise. Obviously, the choice of other normalizations at $z=0$ would lead to a limit profile behaving accordingly in that point.

\subsubsection{Asymptotic behavior for $\gamma \to 0^+$}

As we have underlined at the end of the Introduction, it is well known (see, e.g., \cite{CarKie, Gib}) that \eqref{Mink} reproduces Maxwell dynamics in the limit for $\beta\to +\infty$ (that is, if the maximal field strength is arbitrarily large), and the behavior to be expected is that of the classical theory of Electromagnetism. In terms of \eqref{eqBI}, we thus expect that the critical profiles converge, for $\gamma \to 0^+$, to the critical profile for the linear equation
\begin{equation}\label{lineare}
u_t=u_{xx}+f(u);
\end{equation}
this appears quite natural because we can imagine that the second-order operator in \eqref{eqBI} converges to the second-order derivative for $\gamma \to 0^+$. In order to prove such a convergence, 
we recall that for the second-order linear PDE \eqref{lineare} the picture for traveling fronts can be summarized as follows (see for instance \cite{AroWei, FifMcL}):
\begin{itemize}
\item if $f$ is of type A, then the admissible speeds for \eqref{lineare} form an unbounded interval $[c^*_L, +\infty)$, for a suitable $c^*_L > 0$. 
If $f$ satisfies \eqref{KPP}, then $c_L^*=2\sqrt{f'(0)}$. We denote by $v_L$ the unique front for \eqref{lineare} traveling at the critical speed $c_L^*$ and satisfying \eqref{normalizzazioni};
\item if $f$ is of type B and C, then there exist a unique admissible speed $c^*_L$ and a unique corresponding solution $u(x, t)=v_L(x+c^*_L t)$ of \eqref{lineare} satisfying $v_L(-\infty)=0$, $v_L(+\infty)=1$ and $v_L(0)=\alpha$. 
\end{itemize}
\begin{theorem}\label{modellonuovo}
For every $\gamma > 0$, let $v_\gamma$ be the unique heteroclinic solution of \eqref{eqBI}, connecting $0$ and $1$, which satisfies \eqref{normalizzazioni}. Then, the following hold:
\begin{itemize}
\item[1)] if $f$ is a type A function satisfying \eqref{KPP} (that is, $f$ is a \emph{KPP reaction term}), then 
\begin{equation}\label{convergenze1b}
\lim_{\gamma \to 0^+} c_\gamma^* = c^*_L =2\sqrt{f'(0)}, \qquad  \lim_{\gamma \to 0^+} v_\gamma = v_L;
\end{equation}
\item[2)] if $f$ is of type B or C, it holds
\begin{equation}\label{convergenze1b2}
\lim_{\gamma \to 0^+} c_\gamma^* = c^*_L, \qquad  
\lim_{\gamma \to 0^+} v_\gamma = v_L.
\end{equation}
\end{itemize}
In both cases, the convergence of $v_\gamma$ to $v_L$ is $C^2$ on the whole $\mathbb{R}$.
\end{theorem}
\begin{proof}
1)\, We first notice that \eqref{KPP} and Remark \ref{KPPrem} allow us to read \eqref{Stimaa} as 
$$
\max\left\{2\sqrt{f'(0)}, \sqrt{\sup_{v \in (0, 1]} \gamma^2 \frac{F(v)^2}{v^2} + 2 \frac{F(v)}{v^2}}\right\} \leq c_\gamma^*  \leq 
\sup_{v \in (0, 1]} \left(\gamma f(v) + 2 \sqrt{\frac{f(v)}{v}}\right);
$$
since $\sup_{v \in (0, 1]} f(v)/v = f'(0)$ and $F(v) \leq f'(0) v^2/2$ in view of \eqref{KPP}, we directly infer the first convergence in \eqref{convergenze1b}. As for the convergence of the profiles, the construction of a pointwise a.e. limit for $\{v_\gamma\}_\gamma$ is immediate, thanks to a standard diagonal procedure (see e.g. \cite{Gar, HilKim}) exploiting the boundedness of $\{v_\gamma\}_\gamma$. Hence, there exists a function $\bar{v}$, taking values in $[0, 1]$, for which $v_\gamma \to \bar{v}$ almost everywhere in $\mathbb{R}$. We then take into account the first-order reduction \eqref{notazioney}, here reading as 
\begin{equation}\label{allimite}
\left\{
\begin{array}{l}
\displaystyle y_\gamma'= c_\gamma^* \frac{\sqrt{y_\gamma(2+\gamma^2 y_\gamma)}}{1+\gamma^2 y_\gamma} - f(v) \vspace{0.1cm}\\
y_\gamma(0)=0, \, y_\gamma(1)= 0, \, y_\gamma > 0 \text{ on } (0, 1). 
\end{array}
\right.
\end{equation}
Since $y_{0, f}^-(v)=\int_v^1 f(s) \, ds$ is a supersolution for $(P_{c_\gamma^*, f})^-$ on the whole interval $[0, 1]$, $y_\gamma$ is uniformly bounded in $L^\infty(0, 1)$. 
This has two consequences: on the one hand, thanks to \eqref{allimite}, also $\{y_\gamma'\}_\gamma$ is uniformly bounded in $L^\infty(0, 1)$, and hence $y_\gamma$ converges uniformly to some limit $\bar{y}$; on the other hand, 
by \eqref{costruzionev} we have that also $\{v_\gamma'\}_\gamma$ is bounded in $L^\infty(\mathbb{R})$. 
By the Arzel\`a-Ascoli Theorem, it follows that $v_\gamma \to \bar{v}$ uniformly on compact subsets of $\mathbb{R}$; moreover, since $\{v_\gamma'\}_\gamma$ is bounded and 
\begin{equation}\label{derII}
v_\gamma''(z) = \sqrt{1-\gamma^2 (v_\gamma')^2}^{3} \big(c_\gamma^* v_\gamma'(z) - f(v_\gamma(z))\big),
\end{equation}
also $v_\gamma'$ converges uniformly (necessarily to $\bar{v}'$) on compact subsets of $\mathbb{R}$. Exploiting again \eqref{derII}, $v_\gamma \to \bar{v}$ locally $C^2$. Multiplying \eqref{eqBI} by $\psi \in C_c^\infty(\mathbb{R})$ and integrating on $\mathbb{R}$, we have 
\begin{equation}\label{allimite2}
- \int_{\mathbb{R}} \frac{v_\gamma'(z)}{\sqrt{1-\gamma^2(v_\gamma'(z))^2}} \psi'(z) \, dz +c_\gamma^*  \int_{\mathbb{R}} v_\gamma(z) \psi'(z) \, dz + \int_{\mathbb{R}} f(v_\gamma(z)) \psi(z) \, dz = 0, 
\end{equation}
whence, in view of the previous considerations, passing to the limit for $\gamma \to 0^+$ finally yields 
$$
- \int_{\mathbb{R}} \bar{v}'(z) \psi'(z) \, dz + 2 \sqrt{f'(0)}  \int_{\mathbb{R}} \bar{v}(z) \psi'(z) \, dz + \int_{\mathbb{R}} f(\bar{v}(z)) \psi(z) \, dz = 0.
$$
Consequently, $\bar{v}$ solves the linear equation $v''-c_L^* v'+f(v)=0$ (being $c_L^*=2\sqrt{f'(0)}$) and satisfies $\bar{v}(-\infty)=0$, $\bar{v}(+\infty)=1$, $\bar{v}(0)=1/2$ (this last equality coming from \eqref{normalizzazioni}). 
By the above recalled results on traveling fronts for \eqref{lineare}, it follows that $\bar{v} \equiv v_L$. 
One then deduces the $C^2$-convergence of $v_\gamma$ to $v_L$ on the whole $\mathbb{R}$ from \cite[Lemma 2.4]{Diek}, applied to both $v_\gamma$ and $v_\gamma'$, and from \eqref{derII}. 
\smallbreak
\noindent
2)\,
If instead $f$ is of type B or C, one does not have a priori information about the limit value of $c_\gamma^*$, which we denote by  $\bar{c}$; however, the boundedness of $y_\gamma$ and the convergence of $v_\gamma'$ follow as in the previous argument. Proceeding as for \eqref{allimite2}, passing to the limit for $\gamma \to 0^+$ one has
$$
- \int_{\mathbb{R}} \bar{v}'(z) \psi'(z) \, dz + \bar{c}  \int_{\mathbb{R}} \bar{v}(z) \psi'(z) \, dz + \int_{\mathbb{R}} f(\bar{v}(z)) \psi(z) \, dz = 0;
$$
moreover, $\bar{v}(-\infty)=0$, $\bar{v}(+\infty)=1$ and $\bar{v}(0)=\alpha$ in view of \eqref{normalizzazioni}. We can thus conclude that $\bar{v} \equiv v_L$ similarly as before; moreover, since in this case there exists a unique admissible speed for \eqref{lineare}, necessarily $\bar{c}=c_L^*$, and consequently \eqref{convergenze1b2} is satisfied. Finally, the $C^2$-convergence of $v_\gamma$ follows as in the proof of the previous item.
\end{proof}
\begin{remark}
\textnormal{ 
If $f$ is of type A but does not satisfy the KPP requirement \eqref{KPP}, in principle the quantity $\bar{c}=\lim_{\gamma \to 0^+} c_\gamma^*$, which exists in view of the monotonicity of $c_\gamma^*$, cannot be determined explicitly. 
However, $v_\gamma$ converges to the unique increasing heteroclinic $v$ solving $v''-\bar{c}v+f(v)=0$ such that $v(0)=\mathcal{V}_0$, that is, to the unique front profile for the linear diffusion equation which travels at speed $\bar{c}$ and takes value $1/2$ at $z=0$. If $f'(0)$ exists, then it is worth underlining that this profile necessarily corresponds to the critical one as a consequence of Remark \ref{largeroot}, since the fact that $E(y(v))'(0)$ is the largest root of the equation $x^2-cax+af'(0)=0$ is preserved in the limit (unless $\bar{c}=2\sqrt{f'(0)/a}$, in which case the limit profile is the critical one in view of the previous comments). 
In the general case, instead, whether the limit front profile is critical or not seems less immediate to detect, even if in principle the first alternative should be expected. 
}
\end{remark}
\begin{remark}
\textnormal{In particular, Theorem \ref{modellonuovo} holds in case of balanced reactions (necessarily of type C), where it provides the convergence of the corresponding steady profiles to the unique steady profile for the linear equation \eqref{lineare}; indeed, in this case the only admissible speed for \eqref{lineare} is equal to $0$.}
\end{remark}

\subsubsection{Asymptotic behavior for $\gamma \to +\infty$} 

\begin{figure}[h!]
\center\includegraphics[scale=0.9]{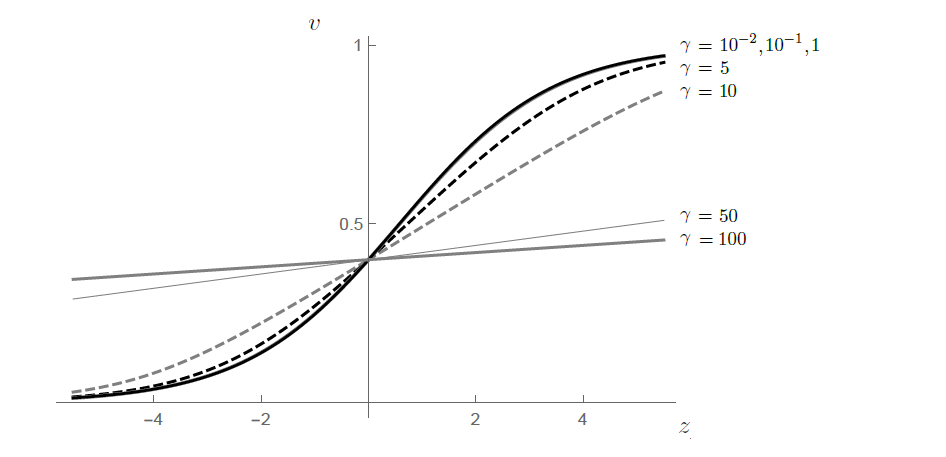}
\caption{Some numerical simulations illustrating Theorems \ref{modellonuovo} and \ref{teorema11eps}, with the positions listed in \eqref{parfigura}}\label{convergenza1}
\end{figure}

We now analyze equation \eqref{eqBI} for $\gamma \to +\infty$, that is, we are interested in the case when the maximal field strength converges to $0$. Here it looks evident that regular profiles have to take small derivatives, in order for the denominator in the second-order term to be well-defined. This naturally leads to the following statement.
\begin{theorem}\label{teorema11eps}
Let $f$ be of type A, B or C. Moreover, let $\{v_\gamma\}_{\gamma}$ be a family of heteroclinic solutions of \eqref{eqBI}, connecting $0$ and $1$, such that $v_\gamma(0) \to v_0 \in (0, 1)$ for $\gamma \to +\infty$.
The following hold:
\begin{itemize}
\item[1.] if $F(1) > 0$, then 
\begin{equation}\label{convergenze1c}
\lim_{\gamma \to +\infty} c_\gamma^* = +\infty, 
\end{equation}
where the convergence occurs monotonically, with order $\gamma$;
\item[2.] if $F(1)=0$ (namely, $f$ is balanced), then 
\begin{equation}\label{convergenze1cbis}
\lim_{\gamma \to +\infty} c_\gamma^* = 0. 
\end{equation}
\end{itemize}
In both cases, it holds 
$$
\lim_{\gamma \to +\infty} v_\gamma = v_0 \quad \text{ in } C^2_{\text{loc}}(\mathbb{R}),
$$
that is, with $C^2$ convergence on compact subsets of $\mathbb{R}$.
\end{theorem}
\begin{proof}
By \eqref{Stimaa} (recall also \eqref{forseutile}), we have that 
$c_\gamma^* \to +\infty \text{ with order } \gamma$ if $F(1) > 0$, so \eqref{convergenze1c} is satisfied; on the contrary, if $F(1)=0$, Remark \ref{bilanciata} ensures that $c_\gamma^* = 0$ for every $\gamma > 0$, hence \eqref{convergenze1cbis} trivially holds. 
As for the front profiles $v_\gamma$, after the construction of a pointwise limit a.e. similar to that in the proof of Theorem \ref{modellonuovo}, the fact that $v_\gamma \in C^2$ for every $\gamma$ directly implies that $\Vert v_\gamma' \Vert_{L^\infty(\mathbb{R})} \leq 1/\gamma$. Consequently, by the Arzel\`a-Ascoli Theorem, $v_\gamma$ converges uniformly to some function $\bar{v}$ on compact sets; moreover, from the equality 
$
v_\gamma(z)-v_\gamma(0)=\int_0^z v_\gamma'(s) \, ds 
$
we deduce that $\bar{v}$ is constant. Thus, $v_\gamma$ converges to the constant value $\bar{v}(0)$, which by assumption coincides with $v_0$. The convergence is proved to be locally $C^2$ in a similar way as for Theorem \ref{modellonuovo}; of course, this time it is not uniform on the whole $\mathbb{R}$, since $v_\gamma(-\infty)=0$ and $v_\gamma(+\infty)=1$ for every $\gamma$.  
\end{proof}
In Figure \ref{convergenza1}, we display some numerical simulations for
\begin{equation}\label{parfigura}
f(s)=s(1-s)(s-0.4), \quad \gamma=10^{-2}, 10^{-1}, 1, 5, 10, 50, 100, 
\end{equation}
where in black we represent the front connecting $0$ and $1$ for equation \eqref{lineare}. Notice that the profiles obtained for $\gamma=10^{-2}$ ($c^* \approx 0.14142$), $\gamma=10^{-1}$ ($c^* \approx 0.14143$), $\gamma=1$ ($c^* \approx 0.14253$) and the critical profile for \eqref{lineare} ($c_L^* \approx 0.14142$) appear indistinguishable; 
indeed, the convergence to $v_L$ for $\gamma \to 0^+$ appears very rapid. On the other hand, the convergence to the constant $\alpha=0.4$ for $\gamma \to +\infty$ has a milder trend. The approximations of the critical speed have been obtained through a shooting procedure which is briefly explained in the short Appendix at the end of the paper; therein, concerning Figure \ref{convergenza1}, we only report the approximated values of $c^*$ from $\gamma=1$ on, since for the considered smaller values of $\gamma$ they differ by quantities smaller than $10^{-4}$. 

\subsection{Behavior of the critical profile in the vanishing (or large) diffusion limit}\label{sez4} 
In this section, we consider the equation
\begin{equation}\label{eq1diff} 
\gamma \left(\frac{v'}{\sqrt{1 - \gamma^2(v')^2}}\right)' - c v' + f(v) =0,
\end{equation}
coinciding with \eqref{TW2} for $a=1/\gamma$, $b=1$, 
for which the diffusion coefficient $\gamma$ is equal to the inverse of the field strength. As for instance in \cite{FolPlaStr}, we can interpret \eqref{eq1diff} as modeling a kind of vanishing (resp., large) diffusion limit for $\gamma \to 0^+$ (resp., $\gamma \to +\infty$), since the velocity $u_x$ is reduced (resp., accelerated) by a factor of $\gamma$; notice that the same effect would be obtained by scaling the spatial variable by a factor of $1/\gamma$. For vanishing diffusion, as we will see, in terms of the behavior of the traveling fronts we will here recover outcomes which are typical of the linear/saturating case, the limit profile being stepwise; this is possible since $\gamma \to 0^+$ and hence the slope of the critical profile may become arbitrarily large. On the other hand, for $\gamma \to +\infty$ the picture is similar to the one in Theorem \ref{teorema11eps}. 

\subsubsection{Asymptotic behavior for $\gamma \to 0^+$}\label{prima}

We first consider \eqref{eq1diff} for $\gamma \to 0^+$, dealing with its critical profile $v_\gamma$ satisfying \eqref{normalizzazioni}. 
Denoting by $H$ the Heaviside function
$\displaystyle
H(z)=\left\{
\begin{array}{ll}
0 & z < 0 \\
1 & z > 0
\end{array}
\right.
$,
we have the following.

\begin{theorem}\label{conv2}
Let $f$ be of type A, B or C. Moreover, for every $\gamma > 0$, let $v_\gamma$ be the unique heteroclinic solution of \eqref{eq1diff}, connecting $0$ and $1$, which satisfies \eqref{normalizzazioni}. 
Then, it holds 
$$
\lim_{\gamma \to 0^+} c_\gamma^* = 0, \qquad \lim_{\gamma \to 0^+} v_\gamma(z) = H(z) \text{ for every } z \neq 0. 
$$
More precisely, the convergence of $c_\gamma^*$ occurs monotonically, with order $\sqrt{\gamma}$; the convergence of $v_\gamma$  is uniform in $\mathbb{R} \setminus \mathcal{I}_0$, $\mathcal{I}_0$ being an arbitrary neighborhood of the origin. 
\end{theorem}
\begin{proof}
In view of \eqref{Stimaa}, \eqref{Stimac} and Proposition \ref{monotono}, it is clear that $c_{\gamma}^* \to 0$ decreasingly and with the order indicated in the statement (recall that $1/a=\gamma$). 
As for the critical profile, we proceed as in the statements of the previous subsection to 
construct a function $\bar{v}$ (taking values in $[0, 1]$)
for which $v_{\gamma} \to \bar{v}$ almost everywhere.  
However, denoted by $y_\gamma$ the associated solution of the first-order reduction \eqref{notazioney}, defined in \eqref{costruzionev}, here the fact that $y_\gamma$ is bounded does not imply that $y'_\gamma$ is bounded as well, since $y_\gamma'=(c_\gamma^*/\gamma) R_{1/\gamma, 1}(y_\gamma) - f(v)$ and $\lim_{\gamma \to 0^+} (c_\gamma^*/\gamma)=+\infty$ (being $c_\gamma^* \sim \sqrt{\gamma}$ for $\gamma \to 0^+$). 
Nevertheless, from \eqref{costruzionev} the boundedness of $y_\gamma$ implies that $\gamma v_\gamma' \to 0$ uniformly for $\gamma \to 0^+$. 
Multiplying \eqref{eq1diff} by $\psi \in C_c^\infty(\mathbb{R})$ and integrating by parts so as to obtain 
\begin{equation}\label{relazII}
- \int_{\mathbb{R}} \frac{\gamma v_{\gamma}'(z)}{\sqrt{1-\gamma^2(v_{\gamma}'(z))^2}} \psi'(z) \, dz +c_{\gamma}^*  \int_{\mathbb{R}} v_{\gamma}(z) \psi'(z) \, dz + \int_{\mathbb{R}} f(v_{\gamma}(z)) \psi(z) \, dz = 0,
\end{equation}
this implies that the first integral in \eqref{relazII} vanishes for $\gamma \to 0^+$; the second integral vanishes as well since $c_\gamma^* \to 0$ for $\gamma \to 0^+$ and $v_\gamma$ is bounded. Consequently,
$$
\lim_{\gamma \to 0^+} \int_{\mathbb{R}} f(v_\gamma(z)) \psi(z) \, dz = 0,
$$
for every $\psi \in C_c^\infty(\mathbb{R})$, implying $f(\bar{v}(z)) = 0$ for almost every $z \in \mathbb{R}$ in view of the boundedness of $f$. The proof now continues differently, according to the shape of the reaction term: 
\smallbreak
\noindent
$\bullet$
If $f$ is of type A, the pointwise convergence follows as in \cite{Gar}: by the monotonicity of $v_\gamma$, together with the fact that $\bar{v}(0)=1/2$, one has that
$$
\bar{v}(z) = 0 \textrm{ for almost every } z < 0, \qquad \bar{v}(z) = 1 \textrm{ for almost every } z > 0, 
$$
and the full pointwise convergence follows from the fact that $z \mapsto v_\gamma(z)$ is monotone for every $\gamma$, together with the comparison theorem for limits. The uniform convergence outside a neighborhood of $0$ follows instead from \cite[Lemma 2.4]{Diek}.
\smallbreak
\noindent
In case $f$ is of type B or C, in principle $v_\gamma$ could instead converge to some other equilibrium, even along a whole interval. In order to show that this cannot occur, we reason as follows: 
\smallbreak
\noindent
$\bullet$
If $f$ is of type B, problem $(P_{c, f})^+$ in \eqref{cauchyf} is explicitly integrable as long as $v \leq \alpha$, yielding
\begin{equation}\label{tipoBesplicita}
y_\gamma(v)=\frac{1}{\gamma}(\sqrt{1+(c_\gamma^*)^2 v^2}-1);
\end{equation}
via \eqref{costruzionev}, together with $v_\gamma(0)=\alpha$, one then gets the implicit relation
\begin{equation}\label{rappresentazione}
\sqrt{1+(c_\gamma^*)^2 v_\gamma(z)^2} +\frac{1}{2} \log \left(\frac{\sqrt{1+(c_\gamma^*)^2 v_\gamma(z)^2}-1}{\sqrt{1+(c_\gamma^*)^2 v_\gamma(z)^2}+1}\right)=\frac{c_\gamma^*}{\gamma} z + k_\gamma c_\gamma^*, 
\end{equation}
with 
$$
k_\gamma c_\gamma^* = \sqrt{1+(c_\gamma^*)^2\alpha^2} +\frac{1}{2} \log \left(\frac{\sqrt{1+(c_\gamma^*)^2\alpha^2} -1}{\sqrt{1+(c_\gamma^*)^2 \alpha^2} +1}\right).
$$
The representation given by \eqref{rappresentazione} is valid until $v_\gamma$ takes value $\alpha$, namely for $z \leq 0$ (recall that $v_\gamma'(z) > 0$ for every $z$). 
We now fix $z \in \mathbb{R}$ and examine the limit of $v_{\gamma}(z)$ as $\gamma \to 0^+$.
\begin{itemize}
\item[-] If $z=z_- < 0$, for $\gamma \to 0^+$ the right-hand side in \eqref{rappresentazione} goes to $-\infty$ with order $1/\sqrt{\gamma}$. Hence, it has necessarily to be $v_\gamma(z_-) \to 0$, otherwise the left-hand side therein would go to $-\infty$ with the slower order $\log((c_\gamma^*)^2)=\log(\gamma)$. It follows that $v_\gamma(z_-) \to 0$, excluding in particular the convergence to other equilibria. 
\item[-] For $z=0$, our normalization condition necessarily implies $v_\gamma(0) \to \alpha$. 
\item[-] Assume now by contradiction that there exists $z=z_+ > 0$ for which $v_\gamma(z_+) \to \alpha$. The fact that $v_\gamma$ is increasing for every $\gamma$ then implies that $\lim_{\gamma \to 0^+} v_\gamma(z) = \alpha$ for every $z \in [0, z_+]$ and also that, fixed $\delta > 0$, taking $\gamma$ sufficiently small it holds $\vert v_\gamma(z)-\alpha \vert \leq \delta$ for every $z \in [0, z_+]$. However, from \eqref{tipoBesplicita} we deduce that for every $v_0 \in (0, \alpha]$ it holds
$$
y_\gamma(v_0) \sim \frac{(c_\gamma^*)^2 v_0^2}{2\gamma} \; \text{ for } \gamma \to 0^+,
$$
whence there exist two constants $\eta, \sigma > 0$, independent of $\gamma$, such that, for every $\gamma > 0$, it holds $y_\gamma(v) > \eta$ for every $v \in [\alpha-\sigma, \alpha]$. Recalling that $f(\alpha)=0$, one can also assume - shrinking $\sigma$, if necessary - that $f(v) < (c_\gamma^*/\gamma) (R_{1/\gamma,1}(\eta/2)/2)$ for $v \in [\alpha, \alpha+\sigma]$; it follows that there exists a positive constant $m$ for which $y_\gamma'(v) > m$ for $v \in [\alpha-\sigma, \alpha+\sigma]$, for every $\gamma > 0$. Therefore, $y_\gamma$ is increasing on $[\alpha-\sigma, \alpha+\sigma]$ and hence $y_\gamma(v) > \eta$ for every $v \in [\alpha-\sigma, \alpha+\sigma]$. 
Recalling \eqref{costruzionev}, we can then find a positive constant $\ell$ such that  $v_\gamma' \geq \ell$ whenever $v_\gamma$ takes values between $\alpha-\sigma$ and $\alpha+\sigma$. 
Hence, 
$\displaystyle
v_\gamma(z_+)-v_\gamma(0)=\int_{0}^{z_+} v_\gamma'(s) \, ds \geq \ell z_+,
$
which is a contradiction passing to the limit for $\gamma \to 0^+$. 
\end{itemize}
Summing up, $\bar{v}$ cannot take value $\alpha$ in points other than $0$ and the fact that $f(\bar{v}(z))=0$ implies $v_\gamma(z) \to 1$ for every $z > 0$. The full pointwise convergence and the uniform convergence then follow as for the type A case.
\smallbreak
\noindent
$\bullet$ If $f$ is of type C, by the previous arguments $v_\gamma$ can converge either to $0$, or to $\alpha$ or to $1$. Here it suffices to show again that $v_\gamma(z) \to \alpha$ if and only if $z=0$. Since $y_\gamma'(v) \geq (c_\gamma^*/\gamma) R_{1/\gamma, 1}(y_\gamma(v))$ for $v \leq \alpha$, one has $y_\gamma(\alpha) \geq \max\Big\{\tfrac{\sqrt{1+(c_\gamma^*)^2 \alpha^2}-1}{\gamma}, -\int_0^\alpha f(s) \, ds\Big\}$. Similarly as before, one can then infer that there exists $\sigma > 0$ such that $v_\gamma' \geq \ell$ whenever $v_\gamma$ takes values between $\alpha-\sigma$ and $\alpha+\sigma$; the above argument for type B functions then works in the same way. The thesis follows. 
\end{proof}
Using the same arguments as in \cite{Gar}, it is possible to see that $v_\gamma' \to \delta_0$ (where $\delta_0$ is the Dirac delta distribution concentrated at $z=0$) in distributional sense; we omit the details for briefness.

\subsubsection{Asymptotic behavior for $\gamma \to +\infty$}

We now deal with equation \eqref{eq1diff} for $\gamma \to +\infty$ (large diffusion limit). Here 
the dynamics is similar to the one in Theorem \ref{teorema11eps}, as we state in 
the following result. 
\begin{theorem}\label{teoremabis}
Let $f$ be of type A, B or C. Moreover, for every $\gamma > 0$, let $v_\gamma$ be the unique heteroclinic solution of \eqref{eq1diff}, connecting $0$ and $1$, satisfying \eqref{normalizzazioni}. The following hold:
\begin{itemize}
\item[1.] if $F(1) > 0$, then 
\begin{equation}\label{convergenze2c}
\lim_{\gamma \to +\infty} c_\gamma^* = +\infty,
\end{equation}
where the convergence occurs monotonically and with order $\gamma$;
\item[2.] if $F(1)=0$ (namely, $f$ is balanced), then 
\begin{equation}\label{convergenze2cbis}
\lim_{\gamma \to +\infty} c_\gamma^* = 0. 
\end{equation}
\end{itemize}
In both cases, it holds
$$
\lim_{\gamma \to +\infty} v_\gamma = \mathcal{V}_0 \quad \text{ in } C^2_{\text{loc}}(\mathbb{R}),
$$
that is, with $C^2$ convergence on compact subsets of $\mathbb{R}$. 
\end{theorem}
\begin{proof}
Since, with the notation in Section \ref{sez2}, it is $b/a=\gamma$, a direct application of Proposition \ref{proponuova} gives that 
$
c_\gamma^* \to +\infty
$
with order $\gamma$ if $F(1) > 0$, so \eqref{convergenze2c} is satisfied; on the other hand, the balanced case provides again  $c_\gamma^*=0$ for every $\gamma$, so that \eqref{convergenze2cbis} is trivial. 
As for the traveling profiles $v_\gamma$, one can proceed as in the proof of Theorem \ref{teorema11eps}, exploiting the fact that $v_\gamma \in C^2$ implies $\Vert v_\gamma' \Vert_{L^\infty(\mathbb{R})} \leq 1/\gamma$, so that $v_\gamma$ converges uniformly to a constant on compact sets. Such a constant is necessarily $\mathcal{V}_0$, the one given by the normalization condition at $z=0$. 
\end{proof}

\begin{remark}\label{conclusivo}
\textnormal{So far, we have actually dealt with cases in which one of the two parameters $a$ and $b$ is equal to $1$. The same outcomes would be obtained if $1$ is replaced by any fixed positive constant $q$, with obvious modifications. Indeed, it would suffice to multiply the entire equation by $q$ to bring the study back to the previous considerations; in this way, obviously $c^*$ would be replaced by $qc^*$ and $f$ by $qf$, but the qualitative features of the results of Subsections \ref{sez3} and \ref{sez4} would continue holding all the same. Referring to Table \ref{tabella} in the Introduction, this is the reason why a generic constant $q$ appears.} 
\end{remark}

\begin{remark}
\textnormal{In the type A case, it is clear that one could also wonder what is the behavior of \emph{noncritical} front profiles, that is, fronts which propagate with speed larger than the critical one. If one fixes the propagation speed $c$, the situations considered in Theorems \ref{teorema11eps} and \ref{teoremabis} have to be excluded, since $c$ ceases to be admissible for sufficiently large $\gamma$. Proceeding as in the proof of Theorem \ref{modellonuovo}, one could instead infer that noncritical profiles traveling at fixed speed $c$ for equation \eqref{eqBI} converge, for $\gamma \to 0^+$, to the noncritical profile $v_L$ of \eqref{lineare} which travels with the same speed $c$. As for \eqref{eq1diff} with $\gamma \to 0^+$, here the behavior of noncritical front profiles traveling at speed $c$ is less trivial: proceeding as for Theorem \ref{conv2}, however, \eqref{relazII} here implies that $v_\gamma$ converges to the linear inviscid front traveling at speed $c$, i.e., the solution of $c\bar{v}'=f(\bar{v})$ satisfying $v(0)=1/2$ (indeed, the second integral therein does not vanish, since $c$ is fixed). This is a typical result for noncritical fronts (see also \cite{Gar, HilKim}).}
\end{remark}

We conclude the present section showing some numerical simulations for equation \eqref{eq1diff}, with
\begin{equation}\label{parfigura2}
f(s)=s(1-s), \quad \gamma=10^{-4}, 10^{-2}, 1, 10, 100, 
\end{equation}
displayed in Figure \ref{convergenza2}. Again, the approximations of the critical speed, reported in the Appendix at the end of the paper, have been obtained through a shooting procedure. 
\begin{figure}[h!]
\center\includegraphics[scale=0.9]{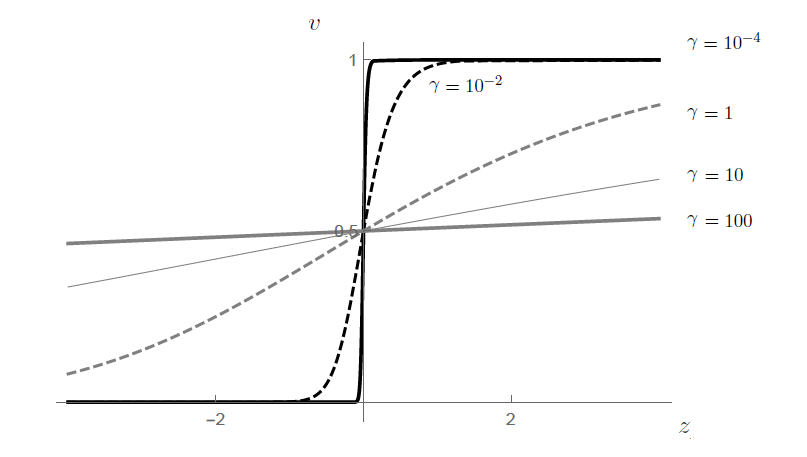}
\caption{Some numerical simulations illustrating Theorems \ref{conv2} and \ref{teoremabis}, with the positions listed in \eqref{parfigura2}}\label{convergenza2}
\end{figure}

\section{New behaviors for the singular perturbation problem} \label{sez5}

In this section, we deal with the $1$-dimensional PDE
$
\displaystyle u_t=\eps \left(\frac{u_x}{\sqrt{1-u_x^2}}\right)_x + f(u),
$
corresponding to \eqref{TW2} with $a=b=1/\eps$.
The critical fronts $v_\eps=v_\eps(z)$ here satisfy the ODE 
\begin{equation}\label{TWdelta}
\eps \left(\frac{v_\eps'}{\sqrt{1-(v_\eps')^2}}\right)' - c_\eps^* v_\eps' + f(v) =0, 
\end{equation}
where $c_\eps^*$ is the critical speed.
Due to the invariance with respect to $z$-translations, we assume as before that every $v_\eps$ satisfies \eqref{normalizzazioni} and we will be interested in the behavior of $v_\eps$ for $\eps \to 0^+$.
As said in the Introduction, we interpret \eqref{TWdelta} as a singular perturbation problem, and indeed we will see that the limit profile never coincides with the inviscid profile $\mathcal{V}_I$ obtained for $\eps=0$. As already mentioned, we could alternatively see \eqref{TWdelta} for $\eps \to 0^+$ as modeling a vanishing diffusion limit, but in view of the terminology adopted in the previous sections we prefer to avoid this interpretation for the sake of clarity. 
\smallbreak
Henceforth,
$y_\eps$ will denote the solution of the first-order reduction \eqref{notazioney}, which here reads as 
\begin{equation}\label{ydeltar}
\left\{
\begin{array}{l}
\displaystyle y_\eps'= c_\eps^* \frac{\sqrt{y_\eps(2\eps+y_\eps)}}{\eps+y_\eps} - f(v) \vspace{0.1cm}\\
y_\eps(0)=0, \, y_\eps(1)= 0, \, y_\eps > 0 \text{ on } (0, 1),
\end{array}
\right.
\end{equation}
and we will set
$$
\bar{c}=\lim_{\eps \to 0^+} c_\eps^*, \quad \bar{v}=\lim_{\eps \to 0^+} v_\eps, \quad \bar{y}=\lim_{\eps \to 0^+} y_\eps,
$$
where we now briefly discuss in which sense such limits should be interpreted.
\smallbreak
Assuming as usual $F(1) > 0$, we first observe that \eqref{Stimaa} and \eqref{Stimac} imply, for every $\eps > 0$, that
\begin{equation}\label{stima2}
F(1) \leq c_\eps^*  \leq \Vert f \Vert_{L^\infty(0, 1)} \quad \text{and} \quad 
F(1) \leq c_\eps^* \leq \frac{1}{\alpha} \sqrt{F(1)-F(\alpha)}.
\end{equation} 
Furthermore, $c_\eps^*$ is increasing with respect to $\eps$, independently of the shape of the reaction term $f$, as a consequence of the fact that $\eps \mapsto \dfrac{\sqrt{y(2\eps+y)}}{\eps + y}$ is decreasing for every $y > 0$. 
Hence, $c_\eps^*$ reaches monotonically its limit $\bar{c}$ for $\eps \to 0^+$ (see also Proposition \ref{monotono}).
Moreover, the inequalities in \eqref{stima2} have the important following consequence.
\begin{proposition}
Let $f$ be of type A, B or C. Then, there exists $\bar{c} > 0$ such that $c_\eps^* \to \bar{c}$ for $\eps \to 0^+$.
\end{proposition}
This means that  
\begin{center}
\textbf{$c_\eps^*$ does not converge to $0$ for $\eps \to 0^+$},
\end{center}
even if the diffusion is arbitrarily slowed down by the parameter $\eps$.
This is a striking difference with respect to both the linear and the saturating diffusion cases (see \cite{Gar}), for which it is well known that $c_\eps^* \to 0$ for $\eps \to 0^+$. 
Incidentally, we notice that one could try to obtain the same conclusion through \eqref{controllomink}, however this could turn slightly more complicated, due to the above discussed lack of homogeneity for \eqref{controllomink} (see also \cite[Remark 3.9]{Gar}). 
\smallbreak
On the other hand, the front profiles $v_\eps$ are equi-Lipschitz continuous and bounded in $C^1$ since $\Vert v_\eps' \Vert_{L^\infty(\mathbb{R})} \leq 1$, 
hence one can construct their limit $\bar{v}$ for $\eps \to 0^+$ similarly as in the previous sections; $\bar{v}$ is a Lipschitz continuous function and $v_\eps \to \bar{v}$ uniformly on compact subsets of $\mathbb{R}$. Finally, for what concerns $y_\eps$, we notice that the backward solution of $y'=-f(v)$ satisfying $y(1)=0$ (given by $y(v)=F(1) - F(v)$) is an upper solution for $y_\eps$. Hence, there exists a constant $Y_\infty > 0$ such that 
\begin{equation}\label{y1}
\Vert y_\eps \Vert_{L^\infty(0, 1)} \leq Y_\infty
\end{equation} 
for every $\eps > 0$; 
since $y \mapsto \dfrac{\sqrt{y(2\eps+y)}}{\eps + y}$ is increasing, the differential equation in \eqref{ydeltar} then yields the existence of a constant $Y_\infty' >0$ for which  
\begin{equation}\label{y2}
\Vert y_\eps' \Vert_{L^\infty(0, 1)} \leq Y_\infty'
\end{equation}
for every $\eps > 0$. 
By the Arzelà-Ascoli Theorem, such a uniform bound on $y'_\eps$ implies that $y_\eps$ converges uniformly to its limit $\bar{y}$ for $\eps \to 0^+$ and $\bar{y}$ is a continuous function, nonnegative on $[0, 1]$. 
\smallbreak
The first difficulty in the study of the asymptotic behavior of $c_\eps^*$ and $v_\eps$ is thus to determine $\bar{c}$. We first briefly comment about the simple case of balanced reaction terms, where everything is explicit, and then we turn to general results for non-balanced reactions. 
 
\subsection{Balanced reaction terms}\label{ibilanciati}

\begin{figure}[h!]
\center\includegraphics[scale=0.75]{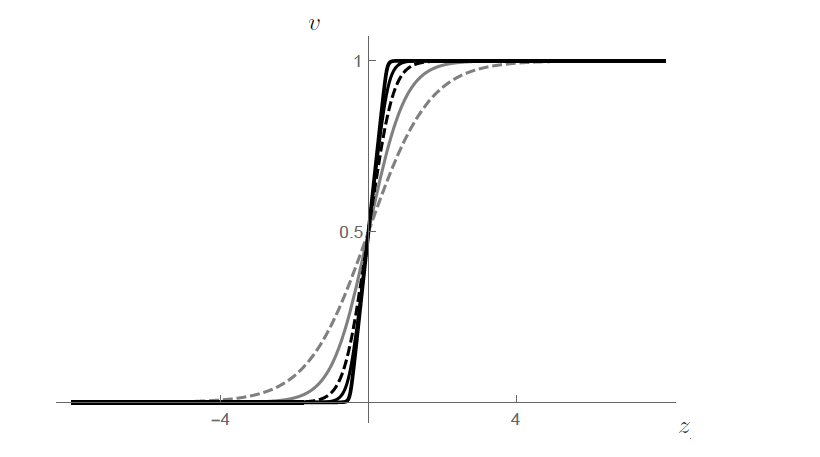}
\caption{Some numerical simulations illustrating the behavior of steady heteroclinic profiles for equation \eqref{TWdelta}, with $f(s)=s(1-s)(s-0.5)$}\label{convergenzab}
\end{figure}

If $f$ is of type C and balanced (that is, $F(1)= 0$), then it is clear that $c_\eps^*=0$ for every $\eps > 0$, so that we are dealing with \emph{steady profiles}. Consequently, $y_\eps(v)=\int_v^1 f(s) \, ds=F(1) -F(v)$ for every $\eps > 0$ and the second differential equation in \eqref{costruzionev} then provides 
$$
v'=\frac{\sqrt{\frac{F(1) -F(v_\eps(z))}{\eps}(2+\frac{F(1) -F(v_\eps(z))}{\eps})}}{1+\frac{F(1) -F(v_\eps(z))}{\eps}}.
$$ 
Therefore, $v_\eps'(z) \to 1$ whenever $v_\eps(z) \not \to 0, 1$ and this necessarily implies that the limit profile $\bar{v}$ is piecewise linear, with slope equal to $1$ whenever taking values different from $0$ and $1$; in view of the normalization condition at $z=0$, one then has $\bar{v}(z)=0$ for $z \in (-\infty, -\alpha]$, $\bar{v}(z)=z+\alpha$ for $z \in [-\alpha, 1-\alpha]$ and $\bar{v}(z)=1$ for $z \in [1-\alpha, +\infty)$. 
In Figure \ref{convergenzab}, for 
\begin{equation}\label{pdelta}
\eps^2=10^{-1} \text{(dashed gray)}, 10^{-2}\text{(gray)}, 10^{-3}\text{(dashed black)}, 10^{-4}\text{(black)}, 10^{-6}\text{(thick black)},
\end{equation}
we represent the steady profiles for $f(s)=s(1-s)(s-0.5)$. 
\\

\subsection{Non-balanced reaction terms}\label{nonbal}

As for non-balanced reaction terms, we start introducing some notation which will be useful in the course. For $\underline{z} \in \mathbb{R}$ and $\underline{v} \in (0, 1)$, we define
$$
\mathcal{V}_L(z; \underline{z}, \underline{v})=\left\{
\begin{array}{ll}
0 & z \in (-\infty, \underline{z}-\underline{v}) \\ 
z-\underline{z}+\underline{v} & z \in [\underline{z}-\underline{v}, 1+\underline{z}-\underline{v}] \\ 
1 & z \in (1+\underline{z}-\underline{v}, +\infty), 
\end{array}
\right.
$$
and, for $f$ of type A, we let $\mathcal{V}_I(z; \underline{z}, \underline{v})$ be the unique (positive) solution of 
$\displaystyle
\left\{
\begin{array}{l}
\bar{c}v'=f(v) \\ 
v(\underline{z})=\underline{v}.
\end{array}
\right.
$
In other words, $\mathcal{V}_L(z; \underline{z}, \underline{v}): \mathbb{R} \to [0, 1]$ is the piecewise linear function which is constant when taking the values $0$ and $1$, has slope $1$ elsewhere and takes value $\underline{v}$ at $z=\underline{z}$, while 
$\mathcal{V}_I(z; \underline{z}, \underline{v})$ is the solution of the autonomous inviscid equation $\bar{c}v'=f(v)$ whose graph passes through the point $(\underline{z}, \underline{v})$. 
Such functions are the ones through which the shape of the limit profile $\bar{v}$ will be constructed. 
Indeed, by the Lipschitz continuity of $\bar{v}$ it makes sense to speak about $\bar{v}'$, at least almost everywhere. Now, 
being $0 \leq v_\eps(z) \leq 1$ and $0 \leq v_\eps'(z) \leq 1$ for every $\eps$ and every $z$ and recalling \eqref{normalizzazioni}, the alternatives which may arise for $\bar{v}$, roughly speaking, are the following:
\begin{itemize}
\item (\emph{non regular scenario}) $\bar{v} \equiv \mathcal{V}_L(\cdot; 0, \mathcal{V}_0)$, namely $\bar{v}$ is fully piecewise linear (with slope equal to $1$ when taking values different from $0$ and $1$);
\item (\emph{regular scenario}) $\bar{v} \equiv \mathcal{V}_I(\cdot; 0, \mathcal{V}_0)$, that is, $\bar{v}$ coincides with the regular inviscid front profile satisfying \eqref{normalizzazioni}
(intuitively, this situation would for instance arise if $v'_\eps$ is bounded away from $1$, so as to make the second-order term in \eqref{TWdelta} disappear for $\eps \to 0^+$); 
\item(\emph{mixed scenario}) 
the limit profile $\bar{v}$ is piecewise linear with slope equal to $1$ on some intervals and coincides with a piece of regular inviscid profile on their complement, in such a way that the junctions between the two different profiles are all continuous; we here avoid writing a formal expression of $\bar{v}$, referring instead to the statement of Theorem \ref{ilteorema} for a case in which the gluing point between the two different profiles is unique.
\end{itemize}
We now want to determine which among the above three possibilities takes place for $\bar{v}$. To begin with, we give a semi-explicit estimate of the limit critical speed for $\eps \to 0^+$ in case the limit profile contains a linear piece with slope $1$. 
\begin{proposition}\label{pezzolineare}
Assume that, for some real numbers $z_0 < z_1$, $v_0 \in [0, 1)$, it holds $\bar{v}(z)=\mathcal{V}_L(z; z_0, v_0)$ on $[z_0, z_1]$
(so that $\bar{v}(z_0)=v_0$) and let $\bar{v}(z_1)=v_1\in (0, 1]$ (hence $\bar{v}(z)=z-z_0+v_0$ on $[z_0, z_1]$, being $z_1-z_0=v_1-v_0$). Then, 
$$
\bar{c}=\frac{F(v_1)-F(v_0)}{v_1-v_0}+\frac{\bar{y}(v_1)-\bar{y}(v_0)}{v_1-v_0}.
$$
\end{proposition}
\begin{proof}
Since $0 \leq v_\eps'(z) \leq 1$ for every $z \in [z_0, z_1]$ and every $\eps > 0$, the family $\{v_\eps'\}_\eps$ is bounded in $L^2(z_0, z_1)$; hence, there exists $w \in L^2(z_0, z_1)$ such that, up to subsequences, $v_\eps' \rightharpoonup w$ in $L^2(z_0, z_1)$ for $\eps \to 0^+$. Extracting a further subsequence, if necessary, $w$ has to coincide almost everywhere with $\bar{v}'$ in the points where $\bar{v}$ is differentiable, and this implies that $w(z)=1$ for almost every $z \in [z_0, z_1]$. Setting $v_\eps^0:=v_\eps(z_0)$, $v_\eps^1:=v_\eps(z_1)$, multiplying now \eqref{TWdelta} by $v'_\eps$ and integrating on $[z_0, z_1]$ one obtains, for every $\eps > 0$,  
\begin{eqnarray}\label{qq1}
c_\eps^* \Vert v_\eps' \Vert^2_{L^2(z_0, z_1)} = c_\eps^* \int_{z_0}^{z_1} v_\eps'(z)^2 \, dz  & = &  \eps\left(\frac{1}{\sqrt{1-v_\eps'(z_1)^2}}-\frac{1}{\sqrt{1-v_\eps'(z_0)^2}}\right)+F(v_1)-F(v_0) \nonumber \\
&  = &  \eps \left(\frac{1}{\sqrt{1-v_\eps'(z_\eps(v_\eps^1))^2}}-\frac{1}{\sqrt{1-v_\eps'(z_\eps(v_\eps^0))^2}}\right) + F(v_1)-F(v_0) \nonumber \\ 
& = & y_\eps(v_\eps^1)-y_\eps(v_\eps^0) + F(v_1)-F(v_0), 
\end{eqnarray}
where $z_\eps$ is the inverse function of $z \mapsto v_\eps(z)$ and we have used \eqref{costruzionev}. 
Since $v_\eps^0 \to v_0$ and $v_\eps^1 \to v_1$, thanks to \eqref{y1} and \eqref{y2} one has that $y_\eps(v_\eps^0) \to \bar{y}(v_0)$ and $y_\eps(v_\eps^1) \to \bar{y}(v_1)$. Hence, in view of the weak semicontinuity of the norm, passing to the inferior limit for $\eps \to 0^+$ on both sides of the above inequality yields
$$
\bar{c} \int_{z_0}^{z_1} w(z)^2 \, dz = \bar{c} (v_1-v_0) \leq \liminf_{\eps \to 0^+} (y_\eps(v_1)-y_\eps(v_0)) + F(v_1)-F(v_0)=\bar{y}(v_1)-\bar{y}(v_0) + F(v_1)-F(v_0).
$$
On the other hand, the fact that $\dfrac{\sqrt{y(2\eps+y)}}{\eps + y} \leq 1$ for every $y > 0$ implies that
\begin{equation}\label{utile1}
y'_\eps(v) \leq c_\eps^* - f(v)
\end{equation}
for every $\eps > 0$ and every $v \in [0, 1]$; 
integrating such an inequality between $v_0$ and $v_1$ and passing to the limit for $\eps \to 0^+$ we obtain 
$
\bar{y}(v_1)-\bar{y}(v_0) \leq \bar{c}(v_1-v_0)- (F(v_1)-F(v_0)),
$
whence we draw the thesis.
\end{proof}
Notice that \eqref{qq1} gives information about the way $v_\eps'$ can converge to $1$, since $\eps/\sqrt{1-v_\eps'(z_\eps(v_\eps))^2}$ is everywhere bounded thanks to \eqref{y1}. 
As a straightforward consequence of Proposition \ref{pezzolineare}, taking into account that $y_\eps(0)=y_\eps(1)=0$ for every $\eps$, we have the following result. 
\begin{corollary}\label{conv1parz}
Let $\bar{v}(z) = \mathcal{V}_L(z; 0, \mathcal{V}_0)$ for every $z \in \mathbb{R}$. 
Then, $\bar{c}=F(1)= \int_0^1 f(s) \, ds$. 
\end{corollary}
In other words, if the non regular scenario takes place for every $z \in \mathbb{R}$ (that is, $v_0=0$ and $v_1=1$), then necessarily the limit critical speed is the lowest possible. 
\smallbreak
On the other hand, we give an obvious necessary condition for the regular scenario to take place on $[0, 1]$. 
\begin{proposition}\label{conv2parz}
Let $\bar{v}(z)=\mathcal{V}_I(z; 0, \mathcal{V}_0)$ for every $z \in \mathbb{R}$. Then, $\bar{c} \geq \Vert f \Vert_{L^\infty(0, 1)}$; if the strict inequality holds, then $\bar{v}'$ is bounded away from $1$, otherwise $\bar{v}$ takes derivative equal to $1$ uniquely when it reaches the value(s) $v$ such that $f(v)=\Vert f \Vert_{L^\infty(0, 1)}$.
\end{proposition}
\begin{proof}
Since $0\leq v_\eps'(z) \leq 1$ for every $z \in \mathbb{R}$, the thesis immediately follows from the fact that $\bar{c} \geq \bar{c} \bar{v}'(z) =f(\bar{v}(z))$ for every $z \in \mathbb{R}$.
\end{proof}
We underline that, due to sign reasons, the situation illustrated in Proposition \ref{conv2parz} may in principle occur only for a type A reaction term, as we will explain at the beginning of Subsection \ref{boc}.
We now divide our study into two different subsections, according to the shape of the reaction term.

\subsection{Reaction terms of type A: one-sided sharpening of the limit profile}

We first deal with type A reactions, for which we introduce the following hypothesis:
\begin{center}
(F) \; $f$ is locally Lipschitz continuous and has a unique local maximum point $v_{\max} \in (0, 1)$.
\end{center} 
The assumption is fulfilled, for instance, in the common cases of Fisher, Huxley or Nagylaki reaction terms; 
notice that it excludes both the existence of subintervals where $f$ is constant and the existence of local minima for $f$ in $[0, 1]$ (and it prescribes a unique change of concavity for $F$, at $v_{\max}$).
For a type A reaction $f$, we have seen in \eqref{stima2} that \eqref{Stimaa} implies
$$
F(1) \leq \bar{c}  \leq \Vert f \Vert_{L^\infty(0, 1)};
$$
the main result of the present section states that, under assumption (F), $\bar{c}$ is strictly included between such bounds and the mixed scenario occurs for the limit profile $\bar{v}$. 
Precisely, we have the following. 
\begin{theorem}\label{ilteorema}
Let $f$ be a type A reaction term fulfilling assumption (F) 
and let $v_+ \in (0, 1)$ be the largest solution of the equation $F(v)=v f(v)$. 
Then, for $\eps \to 0^+$ it holds that $c_\eps^* \to f(v_+)$ and $v_\eps \to \bar{v}$ uniformly in $\mathbb{R}$, where  
$$
\begin{array}{lll}
\bar{v}(z)
 = 
\left\{
\begin{array}{l} 
\mathcal{V}_L(z; 0, \mathcal{V}_0)  \\
\mathcal{V}_I(z; v_+-\mathcal{V}_0, v_+) 
\end{array}
\right.
& 
\begin{array}{l} 
 z\in (-\infty, v_+-\mathcal{V}_0] \\ 
 z \in (v_+-\mathcal{V}_0, +\infty)
\end{array}
&
\text{if } v_+ \geq \mathcal{V}_0 \vspace{0.1cm}\\
\bar{v}(z)
= 
\left\{
\begin{array}{ll} 
\mathcal{V}_L(z; z_+, v_+)\\ 
\mathcal{V}_I(z; 0, \mathcal{V}_0) 
\end{array}
\right.
& 
\begin{array}{l}
 z\in (-\infty, z_+] \\ 
z \in (z_+, +\infty)
\end{array}
& \text{if } v_+ < \mathcal{V}_0,
\end{array}
$$
where, in the latter case, $z_+ < 0$ is defined as the unique real number for which $\mathcal{V}_I(z_+; 0, \mathcal{V}_0)=v_+$.
\end{theorem}
We incidentally notice that the equation $F(v)=vf(v)$ always has at least a positive solution, since $G(v):=F(v)-vf(v)$ is such that $G(0)=0$, $G'(v) < 0$ in a right neighborhood of $0$ and $G(1)>0$. Theorem \ref{ilteorema} thus states that the limit profile $\bar{v}$ is linear - with slope equal to $1$ when strictly positive - up to the value $v_+$, and then solves the inviscid equation $\bar{c}v'=f(v)$. Notice that the gluing between the two profiles takes place in a $C^1$ way, since $\bar{c}=f(v_+)$. 
Hence, a type A reaction term fulfilling (F) has the effect of sharpening \emph{on the left side only} the limit front profile $\bar{v}$ for \eqref{TWdelta}; this asymmetry may appear counterintuitive, especially in case of the Fisher reaction term $f(s)=ms(1-s)$, but it seems hard to lead the behavior of the critical front close to the value $0$ back to its behavior close to $1$ through a simple change of variables, and indeed Theorem \ref{ilteorema} states that this cannot be done. 
\\ 
The proof of Theorem \ref{ilteorema} is given through the following lemmas. 
\begin{lemma}\label{lemma1}
Let $f$ be a type A reaction term fulfilling (F). Then, there exists $v_0 \in (0, 1)$ such that $\bar{y}(v) > 0$ for every $v \in (0, v_0]$; moreover, for every $[\alpha, \beta] \subset (0, v_0)$, $y_\eps \to \bar{y}$ in $C^1([\alpha, \beta])$, where $\bar{y}(v)=\bar{c}v-F(v)$. 
\end{lemma}
\begin{proof}
Following the discussion in Remark \ref{finoaK}, we show that if we suitably choose $v_0 \in (0, 1)$, it is possible to construct a positive lower solution of \eqref{ydeltar} on $[0, v_0]$. 
First, notice that sending $v_0$ to $0$ the estimate \eqref{velocparziale} is fulfilled if
\begin{equation}\label{parzd1}
c \geq \max_{v \in (0, v_0]} f(v)+ 2 \sqrt{\eps\sup_{v \in (0, v_0]}f(v) /v};
\end{equation}
hence, fixed $v_0$ sufficiently small so that $3\max_{v \in (0, v_0]} f(v) < \bar{c}$ and recalling assumption (H), it is possible to find $\eps_0$ such that, for any $0 < \eps < \eps_0$, one has 
$$
\max_{v \in (0, v_0]} f(v)+ 2 \sqrt{\eps\sup_{v \in (0, v_0]}f(v) /v} < 3\max_{v \in (0, v_0]} f(v) < \bar{c},
$$
namely 
$c=3\max_{v \in (0, v_0]} f(v)$ satisfies \eqref{parzd1} and is strictly smaller than $\bar{c}$. With this choice, the requirement \eqref{relazione} will be fulfilled for every $v \in [0, v_0]$ if we take, for instance, $\beta=\max_{v \in (0, v_0]} f(v)$; namely,
$$
y(v)=\eps\left(\sqrt{1+\frac{\beta^2 v^2}{\eps^2}}-1\right)
$$
is a lower solution for $y_\eps$ on the interval $[0, v_0]$, for every $\eps$ sufficiently small. However, fixed $v_1 \in (0, v_0]$, one has 
$$
\eps\left(\sqrt{1+\frac{\beta^2 v_1^2}{\eps^2}}-1\right) \geq \beta v_1-\eps > 0, 
$$
if $\eps$ is sufficiently small.
Hence, in none of the points of the interval $(0, v_0]$ the limit $\bar{y}$ can vanish. Given $[\alpha, \beta] \subset (0, v_0)$, in view of the uniform convergence of $y_\eps$ it follows that $y_\eps$ is uniformly bounded away from $0$ on $[\alpha, \beta]$, and hence 
the right-hand side of the differential equation in \eqref{ydeltar} converges uniformly to $\bar{c}-f(v)$. Since $y_\eps'$ necessarily converges to $\bar{y}'$, we then infer the $C^1$-convergence in the statement.
\end{proof} 
Denoting by $z_{v_0} \in \mathbb{R}$ the point for which $\bar{v}(z_{v_0})=v_0$, Lemma \ref{lemma1} has the important consequence that $\bar{v}(z)=\mathcal{V}_L(z; z_{v_0}, v_0)$ for every $z \in (-\infty, z_{v_0})$, and this excludes the possibility that the regular scenario takes place. 
To see this, notice that for every $z \in (-\infty, z_{v_0})$ such that $\bar{v}(z) > 0$ it holds $y_\eps(v_\eps(z)) \to \bar{y}(\bar{v}(z)) > 0$, thanks to \eqref{y1} and \eqref{y2} 
and recalling that $\bar{v}$ is increasing.  Hence, by \eqref{costruzionev} we obtain
$$
\lim_{\eps \to 0^+} v_\eps'(z)=\lim_{\eps \to 0^+} \frac{\sqrt{y_\eps(v_\eps(z))(2\eps+y_\eps(v_\eps(z)))}}{\eps+y_\eps(v_\eps(z))}
= 1.
$$
Given $z_0, z \in (-\infty, z_{v_0})$ with $z_0 < z$, $\bar{v}(z_0) > 0$, $\bar{v}(z) > 0$, passing to the limit for $\eps \to 0^+$ in the equality $\displaystyle v_\eps(z)-v_\eps(z_0) = \int_{z_0}^z v_\eps'(z) \, dz$  
thanks to the dominated convergence theorem then yields
$
\bar{v}(z)-\bar{v}(z_0)=z-z_0, 
$
as claimed (in particular, $\bar{v}' \equiv 1$ on $(z_0, z)$). 
\smallbreak
Defining now 
\begin{equation}\label{deftildev}
\widetilde{v}:=\sup\{v \in (0, 1] \mid \bar{y} > 0 \text{ on } (0, v)\}, 
\end{equation}
so that $\bar{y}(\widetilde{v})=0$ by the continuity of $\bar{y}$,
we have the following. 
\begin{lemma}\label{lemma1b}
Under the previous assumptions, it holds $\widetilde{v} \in (0, 1)$. Moreover, there exists a unique $z_{\widetilde{v}} \in \mathbb{R}$ such that $\bar{v}(z_{\widetilde{v}})=\widetilde{v}$, and $\lim_{z \to z_{\widetilde{v}}^-} \bar{v}'(z)=1$.
\end{lemma}
\begin{proof}
In view of Lemma \ref{lemma1}, $\widetilde{v} > 0$ and $\bar{y}(v)=\bar{c} v - F(v)$ for every $v \in (0, \widetilde{v}]$ (recall that $\bar{y}$ is continuous). By contradiction, if it were $\widetilde{v}=1$, it would be $\bar{c}=F(1)$ and, in particular, $\bar{v} \equiv \mathcal{V}_L(\cdot; 0, \mathcal{V}_0)$ (compare also with Proposition \ref{pezzolineare}). However, defining $h(v):=F(1)v-F(v)$, we notice that there exists $v_-$ such that $h$ is negative for $v=v_-$, being $h(0)=h(1)=0$ and $h'(v)=F(1)-f(v) > 0$ in a left neighborhood of $v=1$, since $f(1)=0$. This means that $\sup_{v \in (0, 1]} \dfrac{F(v)}{v} > F(1)$, and in view of \eqref{Stimaa} (with $a=b=1/\eps$) this necessarily implies $\bar{c} > F(1)$, yielding a contradiction. Hence, $\widetilde{v} < 1$. 
\\
As for the second part of the statement, were it not the case, then $\bar{v}$ - which is nondecreasing - would be constantly equal to $\widetilde{v}$ on a whole interval $[z_1, z_2]$ ($z_2 > z_1$). We first observe that, similarly as in the proof of Proposition \ref{pezzolineare}, there exists $w \in L^2(z_1, z_2)$ such that $v_\eps \rightharpoonup w$ in $L^2(z_1, z_2)$. 
Multiplying \eqref{TWdelta} by $\psi \in C_c^\infty([z_1, z_2])$ and integrating the first summand by parts on $[z_1, z_2]$ would then yield
$$
- \int_{z_1}^{z_2} \frac{\eps v_{\eps}'(z)}{\sqrt{1-(v_{\eps}'(z))^2}} \psi'(z) \, dz -c_{\eps}^*  \int_{z_1}^{z_2} v_{\eps}'(z) \psi(z) \, dz + \int_{z_1}^{z_2} f(v_{\eps}(z)) \psi(z) \, dz = 0, 
$$
namely, using \eqref{costruzionev},
\begin{equation}\label{distr1}
- \int_{z_1}^{z_2} (y_\eps(v_\eps(z))+\eps) \psi'(z) \, dz -c_{\eps}^*  \int_{z_1}^{z_2} v_{\eps}'(z) \psi(z) \, dz + \int_{z_1}^{z_2} f(v_{\eps}(z)) \psi(z) \, dz = 0.
\end{equation}
For every $z \in [z_1, z_2]$, one has $y_\eps(v_\eps(z)) \to \bar{y}(\widetilde{v})=0$ for $\eps \to 0^+$ and $v_\eps(z) \to \widetilde{v}$, where the convergences are both uniform (the first one thanks to \eqref{y1} and \eqref{y2}); hence, \eqref{distr1} would produce, for $\eps \to 0^+$, 
$\int_{z_1}^{z_2} (-\bar{c}w(z)+f(\widetilde{v})) \psi(z)\, dz =0$. By the arbitrariness of $\psi$, one would deduce that $w(z)=f(\widetilde{v})/\bar{c}$ for every $z \in [z_1, z_2]$, namely $w$ coincides with a positive constant. However, since $\bar{v} \equiv \widetilde{v}$ on $[z_1, z_2]$, passing to the limit for $\eps \to 0^+$ this contradicts the equality $v_\eps(z_2)-v_\eps(z_1) = \int_{z_1}^{z_2} v_\eps'(z) \,dz$.
The limit profile $\bar{v}$ is thus piecewise linear with slope $1$ in a suitable left neighborhood of $z_{\widetilde{v}}$ (i.e., when taking values in $(0, \widetilde{v})$), so that  
$\lim_{z \to z_{\widetilde{v}}^-} \bar{v}'(z)=1$. 
\end{proof}
So far, we have proved that neither the regular nor the non regular scenario occurs. Hence, a mixed picture necessarily occurs, and the remainder of the section is devoted to showing that $\bar{v}$ coincides with the profile in the statement of Theorem \ref{ilteorema}. First, we determine the value of $\bar{c}$. 
\begin{lemma}\label{lemma2}
Let the assumptions of Lemma \ref{lemma1} be fulfilled and let $v_m \in (0, \widetilde{v})$ be such that $\bar{y}(v_m)=M:=\max_{v \in [0, \widetilde{v}]} \bar{y}(v)$. Then, 
\begin{equation}\label{valorebarc}
\bar{c}=\frac{F(\widetilde{v})}{\widetilde{v}}=f(v_m)=f(\widetilde{v}). 
\end{equation}
\end{lemma}
\begin{proof}
The first equality is a straight consequence of \eqref{deftildev}, since $0=\bar{y}(\widetilde{v})=\bar{c}\widetilde{v}-F(\widetilde{v})$. Also the second equality is immediate: by the $C^1$-convergence of $y_\eps$ to $\bar{y}$ on a compact neighborhood of $v_m$ and by \eqref{ydeltar}, for $\eps \to 0^+$ one has 
\begin{equation}\label{interm1}
0 \leftarrow \displaystyle y_\eps'(v_m)= c_\eps^* \frac{\sqrt{y_\eps(v_m)(2\eps+y_\eps(v_m))}}{\eps+y_\eps(v_m)} - f(v_m) \to \bar{c} - f(v_m),
\end{equation}
since $y_\eps(v_m) \to M > 0$. As for the last equality, notice first that, fixed $v_1 \in (\widetilde{v}, 1)$, integrating both sides of \eqref{utile1} between $\widetilde{v}$ and $v_1$ yields 
$
y_\eps(v_1)- y_\eps(\widetilde{v}) \leq c_\eps^*(v_1-\widetilde{v}) - (F(v_1)-F(\widetilde{v})),
$ 
which passing to the limit for $\eps \to 0^+$ produces (recalling that $\bar{y}(\widetilde{v})=0$)
$$
\bar{c} \geq \frac{F(v_1)-F(\widetilde{v})}{v_1-\widetilde{v}};
$$
passing now to the limit for $v_1 \to \widetilde{v}^+$, we obtain the inequality $\bar{c} \geq f(\widetilde{v})$.
As for the reversed inequality, we preliminarily observe that $\bar{y}$ is necessarily decreasing for $v \in (v_m, \widetilde{v})$. Indeed, were this not true, recalling that $\bar{y} \in C^1((0, \widetilde{v}))$ one would have the existence of at least two further critical points of $\bar{y}$ to the right of $v_m$, a local minimum $v_m'$ and a local maximum $v_m''$ (recall that $\bar{y}(\widetilde{v})=0$), where $\bar{y}$ takes strictly positive values (for the local minimum, this holds in view of the definition of $\widetilde{v}$). However, by arguing as for \eqref{interm1}, one would then find that $\bar{c}=f(v_m')=f(v_m'')$, so that the equation $f(v)=\bar{c}$ would possess at least the three distinct solutions $v_m, v_m', v_m''$, all belonging to the interval $(0, \widetilde{v})$. This would contradict assumption (F). Hence, recalling the expression of $\bar{y}$ in the interval $(0, \widetilde{v})$, for every $v \in (v_m, \widetilde{v})$ one has
\begin{equation}\label{osserv}
0 \geq \bar{y}'(v) = \bar{c} - f(v); 
\end{equation}
passing to the limit for $v \to \widetilde{v}^-$, one obtains $\bar{c} \leq f(\widetilde{v})$. The last equality in \eqref{valorebarc} follows and the proof is completed.
\end{proof}
Observe that assumption (F) has been crucial in order to prove the last equality in \eqref{valorebarc}. Moreover, since $v_m < \widetilde{v}$, from the previous proof we can desume that $\widetilde{v}=v_+$, in particular $\widetilde{v}$ lies to the right of the maximum point $v_{\max}$ for $f$ (or, equivalently, of the point where $F$ changes its concavity). Hence, $f$ is \emph{strictly decreasing to the right of $\widetilde{v}$}. 
Incidentally,
from \eqref{osserv} we see in fact that $\bar{c} \leq f(v)$ for every $v \in (v_m, \widetilde{v})$, in accord with the fact that $f$ takes values larger than $f(\widetilde{v})$ in this interval. \\
We are now ready to complete the proof of Theorem \ref{ilteorema}. 
\begin{proof}[Proof of Theorem \ref{ilteorema}.]
For any fixed $\sigma > 0$, we focus on the behavior of $y_\eps$ on $[\widetilde{v}+\sigma, 1]$. Let $M_{\eps, \sigma}:=\max_{v \in [\widetilde{v}+\sigma, 1]} y_\eps(v)$ and assume by contradiction that $\eps/M_{\eps, \sigma} \to 0$. Denoting by $v_{\eps, \sigma} \in  [\widetilde{v}+\sigma, 1)$ the point for which $y_\eps(v_{\eps, \sigma} )=M_{\eps, \sigma}$ (trivially $v_{\eps, \sigma}  < 1$, since $y_\eps(1)=0$), up to subsequences one will have $v_{\eps, \sigma}  \to v^* \in [\widetilde{v}+\sigma, 1]$ and $y_\eps(v_{\eps, \sigma} )=M_{\eps, \sigma} \to \bar{y}(v^*)$, thanks to \eqref{y1} and \eqref{y2}. Hence, it will be 
$$
0 \geq y_\eps'(v_{\eps, \sigma} )=c_\eps^* \frac{\sqrt{y_\eps(v_{\eps, \sigma} )(2\eps+y_\eps(v_{\eps, \sigma} ))}}{\eps+y_\eps(v_{\eps, \sigma} )} - f(v_{\eps, \sigma} )=\displaystyle c_\eps^* \frac{\sqrt{2\frac{\eps}{M_{\eps, \sigma}}+1}}{\frac{\eps}{M_{\eps, \sigma}}+1} - f(v_{\eps, \sigma} ) \to \bar{c}-f(v^*), 
$$
implying $f(\widetilde{v})=\bar{c} \leq f(v^*)$ (recall Lemma \ref{lemma2}), which is a contradiction since $f$ is strictly decreasing to the right of $\widetilde{v}$, as observed above. Consequently, on the one hand, for every $\sigma > 0$ it holds $M_{\eps, \sigma} \to 0$ for $\eps \to 0^+$, so that $\bar{y}(v)=0$ for every $v \in [\widetilde{v}, 1]$, recalling the definition of $\widetilde{v}$; on the other hand, for every $\sigma > 0$, $M_{\eps, \sigma} \to 0$ exactly with order $\eps$,
since if it were $M_{\eps, \sigma}/\eps \to 0$ 
for $\eps \to 0^+$, then \eqref{ydeltar} would imply $y_\eps' \to -f$ uniformly on $[\widetilde{v}+\sigma, 1]$, a contradiction with the fact that $y_\eps \to 0$ uniformly on such an interval. 
Let now $\sigma > 0$ be fixed and let $[z_\sigma, +\infty)=\bar{v}^{-1}([\widetilde{v}+\sigma, 1])$. We claim that there exists $0 < K_\sigma < 1$ such that $\Vert v_\eps' \Vert_{L^\infty(z_\sigma, +\infty)} \leq K_\sigma$ for every sufficiently small $\eps$. Indeed, by the above computations there exists $\hat{M}_\sigma > 0$ such that $M_{\eps, \sigma}/\eps \leq \hat{M}_\sigma$ for every $\eps > 0$, so that for every $z \geq z_\sigma$ it holds
\begin{equation}\label{boundv'}
v_\eps'(z)=\frac{\sqrt{y_\eps(v_\eps(z))(2\eps+y_\eps(v_\eps(z)))}}{\eps+y_\eps(v_\eps(z))} \leq  \frac{\sqrt{\frac{y_\eps(v_\eps(z))}{\eps}\Big(2+\frac{y_\eps(v_\eps(z))}{\eps}\Big)}}{1+\frac{y_\eps(v_\eps(z))}{\eps}}\leq \frac{\sqrt{\hat{M}_{\sigma/2}(2+2\hat{M}_{\sigma/2})}}{1+\hat{M}_{\sigma/2}}=:K_\sigma
\end{equation}
(the use of $\hat{M}_{\sigma/2}$ allows one to control the term $v_\eps(z)$, which converges to $\bar{v}(z)$, on $[z_\sigma, +\infty)$, for $\eps$ sufficiently small).
It follows that $v_\eps'$ is bounded away from $1$ on $[z_\sigma, +\infty)$. Writing \eqref{TWdelta} in distributional form on $[z_\sigma, +\infty)$ as 
$$
- \int_{z_\sigma}^{+\infty} \frac{\eps v_{\eps}'(z)}{\sqrt{1-(v_{\eps}'(z))^2}} \psi'(z) \, dz +c_{\eps}^*  \int_{z_\sigma}^{+\infty} v_{\eps}(z) \psi'(z) \, dz + \int_{z_\sigma}^{+\infty} f(v_{\eps}(z)) \psi(z) \, dz = 0
$$
($\psi \in C_c^\infty([z_\sigma, +\infty))$), this fact allows to infer that for $\eps \to 0^+$ the first summand vanishes and hence $\bar{v}$ solves the inviscid equation $\bar{c} \bar{v}'(z)=f(\bar{v}(z))$ for every $z \geq z_\sigma$; in particular, $\bar{v} \in C^1([z_\sigma, +\infty))$. 
Repeating the argument for every $v > \widetilde{v}$, one actually has that $\bar{c} \bar{v}'(z)=f(\bar{v}(z))$ for every $z > z_{\widetilde{v}}$ (indeed, since  $\bar{c} \bar{v}'(z)=f(\bar{v}(z))$ for every $z > z_{\widetilde{v}}$, it holds that $\bar{v}'(z) > 0$ for every $z > z_{\widetilde{v}}$, hence $\bar{v}$ is strictly increasing therein). Now, by the continuity of $\bar{v}$, one has 
$$
\lim_{z \to z_{\widetilde{v}}^+} \bar{c} \bar{v}'(z) = \lim_{z \to z_{\widetilde{v}}^+}f(\bar{v}(z)) = f(\widetilde{v});
$$
together with Lemma \ref{lemma1b}, this
implies $\bar{v}'(z_{\widetilde{v}})=1$. The proof is completed.  
\end{proof}
\begin{remark}
\textnormal{We notice that \cite[Lemma 2.4]{Diek} ensures that $v_\eps \to \bar{v}$ uniformly on the whole $\mathbb{R}$. Moreover, recalling also \eqref{boundv'}, such a convergence is $C^1$ on each interval $[z, +\infty)$ with $z > z_{\widetilde{v}}$, as well as on each interval $[z_-, z]$, with $z < z_{\widetilde{v}}$, on which $\bar{v} > 0$.  }
\end{remark}

\begin{figure}[h!]
\includegraphics[scale=0.68]{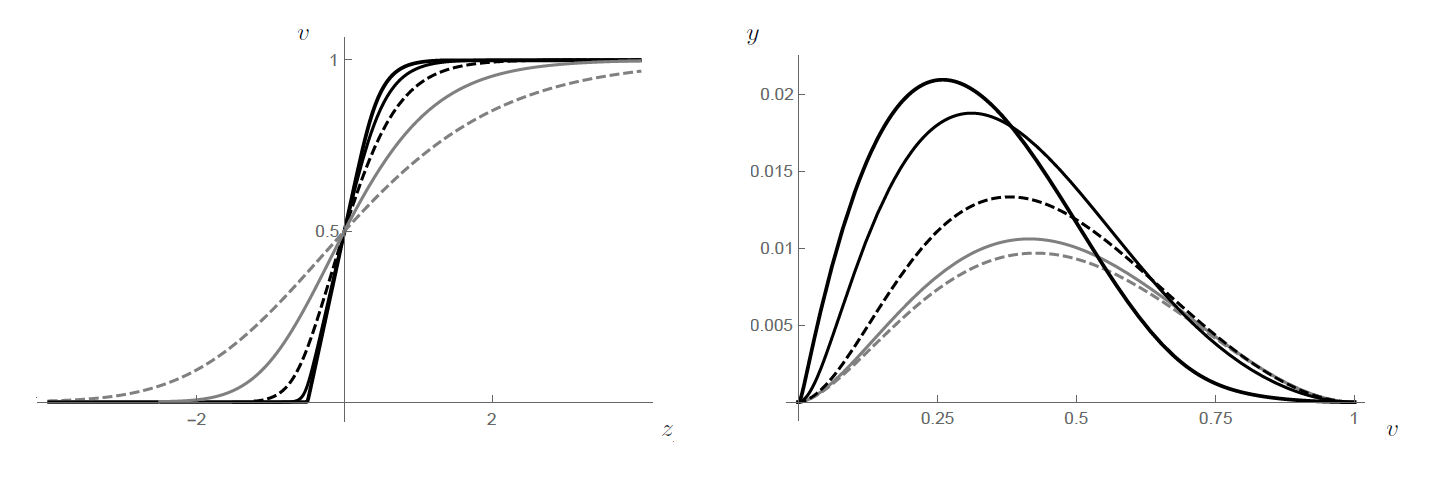}
\caption{Some numerical simulations illustrating the behavior of the critical fronts for equation \eqref{TWdelta} with $f(s)=s(1-s)$. On the left we draw the second-order profile $v$, while on the right we display the corresponding first-order profile $y$ (recall \eqref{notazioney})}\label{fisher}
\end{figure}

\begin{figure}[h!]
\includegraphics[scale=0.68]{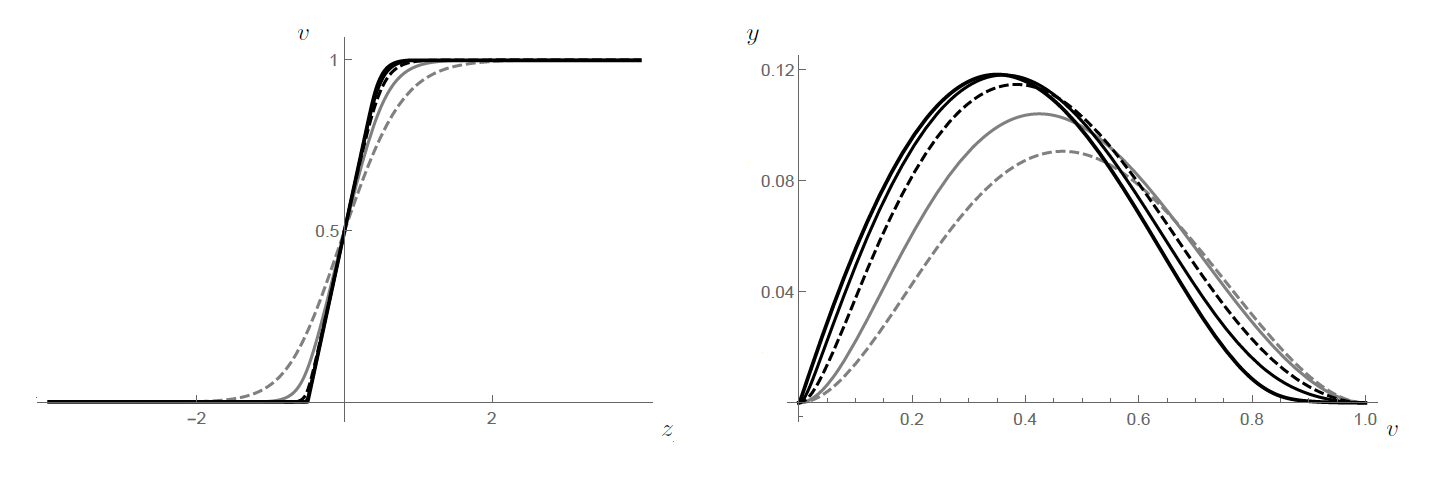}
\caption{Some numerical simulations illustrating the behavior of the critical fronts for equation \eqref{TWdelta} with $f(s)=s(1-s)(1+5s)$. On the left we draw the second-order profile $v$, while on the right we display the corresponding first-order profile $y$ (recall \eqref{notazioney})}\label{nagy}
\end{figure}

We show some numerical simulations for the reaction terms $f(v)=v(1-v)$ (Figure \ref{fisher}) and $f(v)=v(1-v)(1+5v)$ (Figure \ref{nagy}), for which the critical speed has been determined numerically, by a shooting method. We depict the fronts $v_\eps$ and the corresponding first-order profile $y_\eps$ for $\eps$ as in \eqref{pdelta}; the numerical approximations of the critical speed are reported at the end of the paper.
\\
For both the situations displayed in the figure, the shape of the limit profiles predicted by Theorem \ref{ilteorema} arises quite evidently. As we know in view of the previous discussion, the value $v_+$ taken by the limit profile in correspondence of the gluing has to be on the right of the 
maximum of $f$ ($v_{\max}=1/2$ in Figure \ref{fisher} and $v_{\max}=4/15+\sqrt{31}/15 \approx 0.638$ in Figure \ref{nagy}). Computing the solution $v_+$ of the equation $F(v_+)/v_+=f(v_+)$, we obtain $v_+=0.75$ for the Fisher equation and $v_+ \approx 0.865$ for the Nagylaki reaction term; when taking values greater than $v_+$, $\bar{v}$ coincides with $\mathcal{V}_I(\cdot; v_+-\mathcal{V}_0, v_+)$. The right pictures in Figures \ref{fisher} and \ref{nagy} are in accord with this computation, since we see that, starting from $v_+$ on, the first-order profile $y_\eps$ converges to $0$ uniformly, which must be the case (with order $\eps$) for $\bar{v}$ to coincide with an inviscid piece of front profile, in view of the previous discussion. 
Furthermore, the numerical simulations also confirm quite neatly that $\bar{y}$ reaches its maximum in the point $v_m$ such that $f(v_m)=f(v_+)$, as predicted by \eqref{valorebarc}.
Finally, in the $z$-variable, the gluing between the two different profiles here takes place at the point $v_+-1/2$, as the left pictures in the figures seem to confirm. 
\\
Summing up, we have proved that, quite surprisingly, the limit profile sharpens on one side only, being smooth near the other equilibrium. This is the case also for reaction terms of type B, as we will see in the next subsection.

\subsection{Reaction terms of type B}

\begin{figure}[h!]
\includegraphics[scale=0.68]{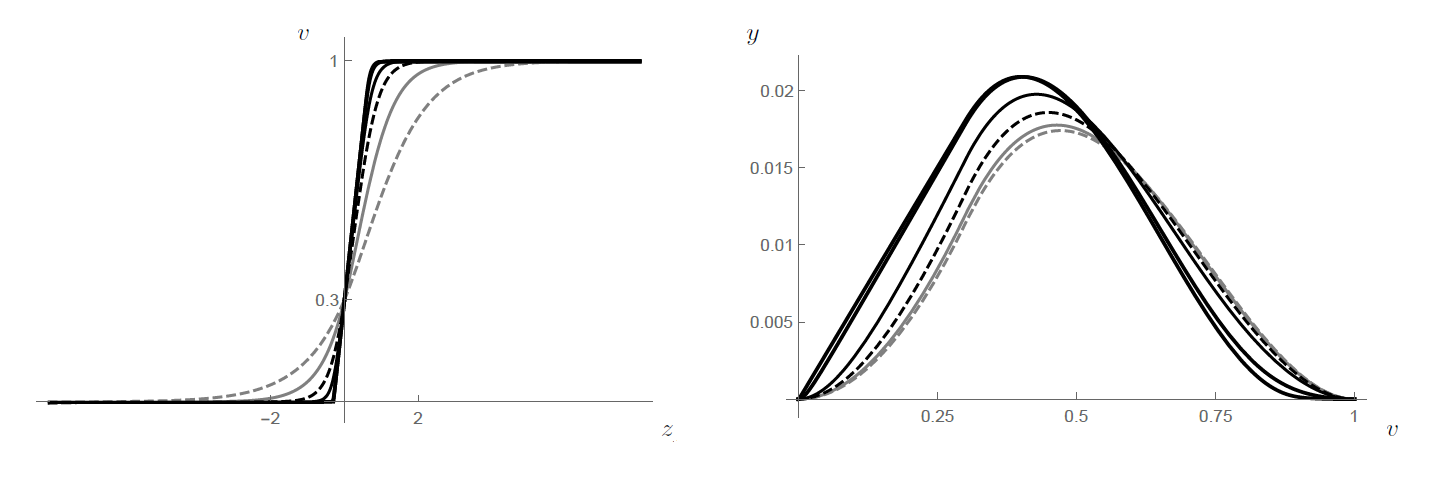}
\caption{Some numerical simulations illustrating the behavior of the critical fronts for equation \eqref{TWdelta} with $f(s)=0$ for $s < 0.3$, $f(s)=(s-0.3)(1-s)$ for $s \geq 0.3$. On the left we draw the second-order profile $v$, while on the right we display the corresponding first-order profile $y$ (recall \eqref{notazioney})}\label{tipoB1}
\end{figure}

As for a type B reaction term $f$, for the sake of simplicity and in order to exploit some of the arguments of the previous section we will assume that $f\big\vert_{[\alpha, 1]}$ satisfies assumption (F). With this preliminary, we first observe that the fact that $\bar{c} > 0$, together with the sign of $f$, implies that in the limit for $\eps \to 0^+$ the first summand in \eqref{TWdelta} cannot go to $0$. This implies that $\bar{v}$ is necessarily piecewise linear with slope $1$ as long as it takes values in $(0, \alpha]$. On the other hand, one can 
repeat the argument of the previous section to infer that $\bar{c}$ cannot be equal to $F(1)$. Therefore, for type B reaction terms the following theorem holds.
\begin{theorem}\label{ilteorematipoB}
Let $f$ be a type B reaction term fulfilling assumption (H) and such that $f\big\vert_{[\alpha, 1]}$ satisfies (F). 
Then, for $\eps \to 0^+$ it holds that
$c_\eps^* \to f(v_+)$, where $v_+ \in (0, 1)$ is the largest solution of the equation $F(v)=v f(v)$, and  
$v_\eps \to \bar{v}$ uniformly in $\mathbb{R}$, where 
$$
\bar{v}(z)
 = 
\left\{
\begin{array}{ll} 
\mathcal{V}_L(z; 0, \mathcal{V}_0)  &  z\in (-\infty, v_+-\mathcal{V}_0] \\
\mathcal{V}_I(z; v_+-\mathcal{V}_0, v_+) & z \in (v_+-\mathcal{V}_0, +\infty),
\end{array}
\right.
$$
(recalling that in this case $\mathcal{V}_0=\alpha$). 
\end{theorem}
In Figure \ref{tipoB1}, we display a numerical experiment illustrating the statement of Theorem \ref{ilteorematipoB}. We see again that the first-order profile $\bar{y}$ has the expected shape, in that for values close to $1$, $y_\eps \to 0$ and indeed $\bar{v}$ is smooth, whereas near $0$ the first-order profile remains positive and in fact $\bar{v}$ is piecewise linear.

\subsection{Reaction terms of type C}\label{boc}

As for type C reaction terms, the same sign argument as for type B functions implies that $\bar{v}$ is necessarily piecewise linear as long as it takes values in $(0, \alpha]$. We can actually push this argument a little bit beyond, as the following proposition states; to this end, we recall that \eqref{normalizzazioni} here prescribes $v_\eps(0)=\alpha$. 
\begin{proposition}
Let $f$ be a  type C reaction term and let $v_* \in (\mathcal{V}_0, 1)$ be defined by the equality $F(v_*)=0$.
Then, for $\eps \to 0^+$ it holds that $v_\eps \to \mathcal{V}_L(\cdot; 0, \mathcal{V}_0)$ uniformly in $(-\infty, v_*-\mathcal{V}_0)$ (recalling that in this case $\mathcal{V}_0=\alpha$).
\end{proposition}
\begin{proof}
Thanks to the position $v_\eps(0)=\mathcal{V}_0=\alpha$, for every $\eps > 0$ it holds $f(v_\eps(z)) < 0$ for every $z < 0$; in view of the sign of the summands in \eqref{TWdelta}, we can thus infer that, for every $z < 0$, either $v_\eps(z) \to 0$ or $v_\eps'(z) \to 1$. Since $v_\eps$ is increasing for every $\eps$, the limit function $\bar{v}$ is continuous and $\bar{v}(0) = \alpha$ due to our normalization condition, necessarily $v_\eps(z) \to \mathcal{V}_L(z; 0, \mathcal{V}_0)$ for every $z \leq 0$. 
Such a convergence can actually be prolonged on a suitable right neighborhood of $z=0$: 
indeed, for every $\eps > 0$ the solution $\hat{y}$ of
$\displaystyle 
\left\{
\begin{array}{l} 
y'=-f(v) \\ 
y(0)=0 
\end{array}
\right.$ 
is a positive strict subsolution for $y_\eps$ on $(0, v_*)$, hence $y_\eps(v) \geq \hat{y}(v) > 0$ for every $v \in (0, v_*)$. For $v < v_*$ one then has $\bar{y}(v) > 0$, implying that $v_\eps'(z)$ converges to $1$ for every $z$ such that $\bar{v}(z) < v_*$. Summing up, $v_\eps(z) \to \mathcal{V}_L(z; 0, \mathcal{V}_0)$ for every $z \in (-\infty, v_*-\mathcal{V}_0)$, whence the thesis.
\end{proof}
\begin{center}
\begin{figure}[h!]
\includegraphics[scale=0.65]{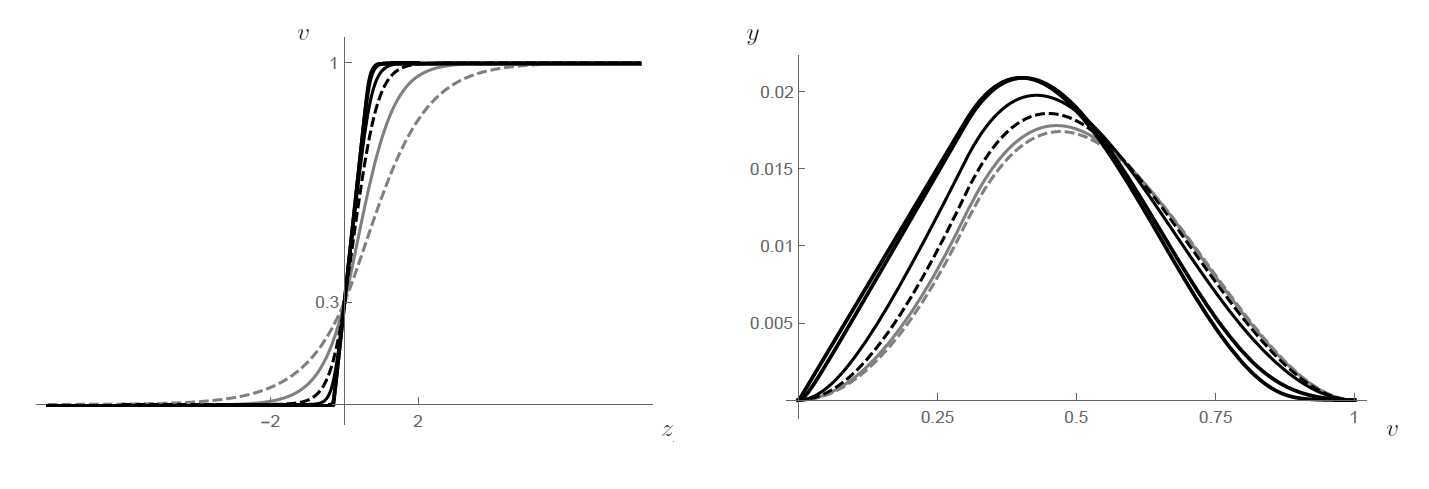}
\caption{Some numerical simulations illustrating the behavior of the critical fronts for equation \eqref{TWdelta}, with $f(s)=10^{-4}s(1-s)(s-0.3)$ for $s < 0.3$ and $f(s)=(s-0.3)(1-s)$ for $0.3 \leq s < 1$. On the left we draw the second-order profile $v$, while on the right we display the corresponding first-order profile $y$ (recall \eqref{notazioney})}\label{convergenzacsm}
\end{figure}
\end{center}

\begin{center}
\begin{figure}[h!]
\includegraphics[scale=0.65]{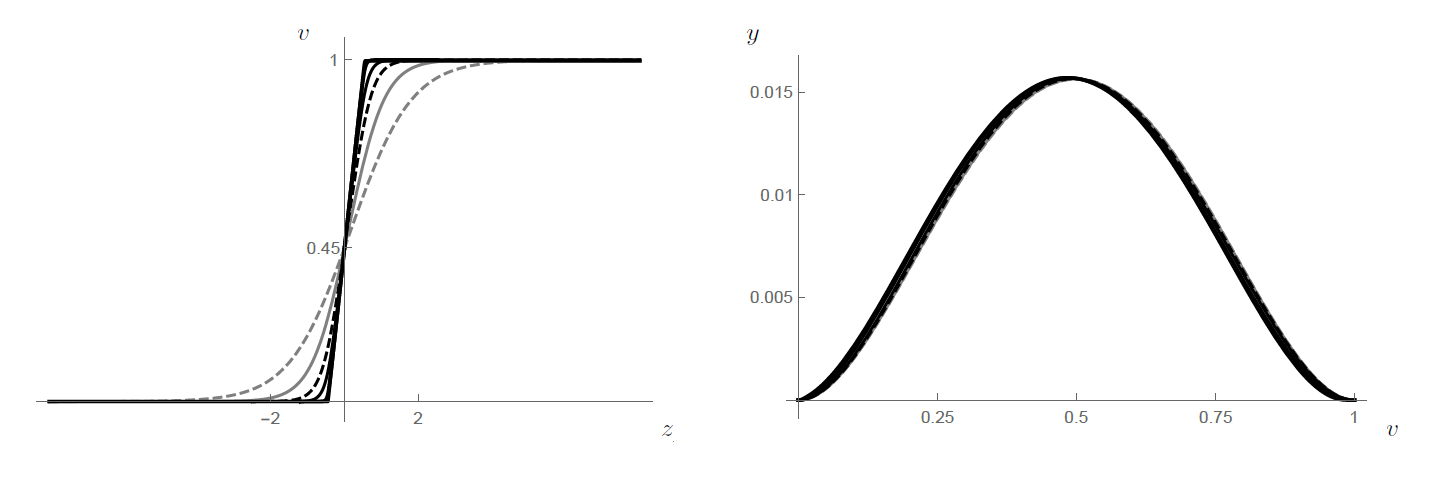}
\caption{Some numerical simulations illustrating the behavior of the critical fronts for equation \eqref{TWdelta}, with $f(s)=s(1-s)(s-0.45)$. On the left we draw the second-order profile $v$, while on the right we display the corresponding first-order profile $y$ (recall \eqref{notazioney})}\label{convergenzacpl}
\end{figure}
\end{center}

In the previous proof, we do not have precise estimates on $v_*$, this being the reason why the statement ensures the convergence to the piecewise linear profile only until $z=v_*-\mathcal{V}_0$; notice moreover that $v_*-\mathcal{V}_0=0$ in the type B case. The possibilities for $v > v_*-\mathcal{V}_0$ are now two: either the limit profile continues smoothly from some point on, or it remains equal to $\mathcal{V}_L(z; 0, \mathcal{V}_0)$ until taking value $1$. Likely, the alternative which takes place depends on the reaction term $f$. On the one hand, indeed, by a continuity argument it appears more than reasonable that if the negative part of $f$ is very small, then we will have a picture similar to the one for type B reactions, as Figure \ref{convergenzacsm} seems to suggest. On the other hand, if $f$ is close to a balanced reaction term, it may probably be expected that the limit profile does not differ much from the one for balanced reaction terms. Intuitively, if the negative part of $f$ is sufficiently large, in the first-order reduction the lower solution from the left $y(v)=-\int_0^v f(s) \, ds$ may in fact reach an arbitrarily large value at $v=\alpha$, so that it may be that $y_\eps$ does not converge to $0$ in any neighborhood of $1$. This would imply that $\bar{v}$ is everywhere piecewise linear, with slope equal to $1$ until taking value $1$, and $c_\eps^* \to F(1)$ (as seems the case in Figure \ref{convergenzacpl}). We thus conjecture that \emph{both the behaviors are possible}, in dependence on the reaction; nevertheless, it does not seem immediate to reach a precise statement in this respect, so we just limit our discussion to the pictures in Figures \ref{convergenzacsm} and \ref{convergenzacpl}. 

\section{Final remarks}\label{sez6} 

\begin{center}
\begin{figure}[h!]
\includegraphics[scale=0.46]{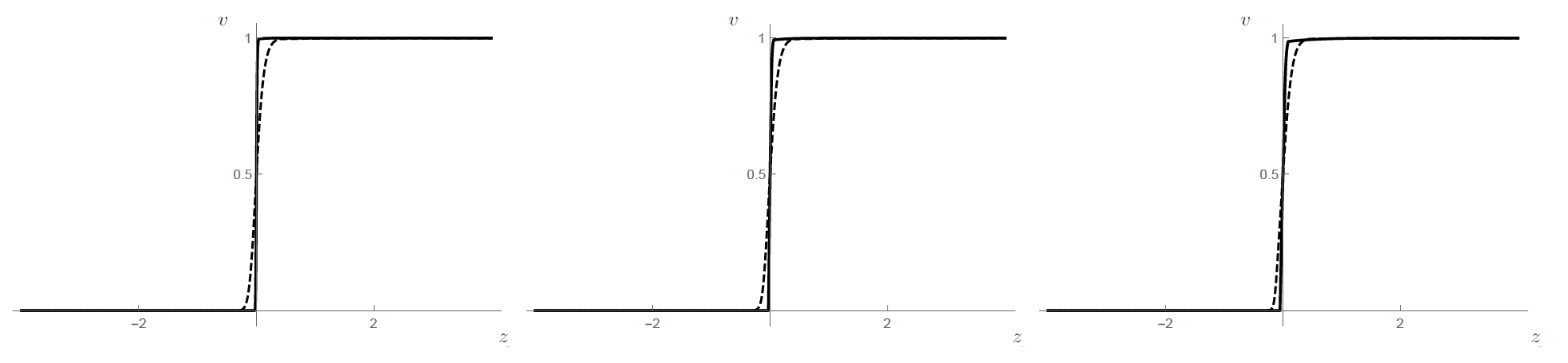}
\caption{Some numerical simulations illustrating the behavior of the critical fronts for equation \eqref{general+}-$(i)$, with $f(s)=s(1-s)$}\label{convergenzalast}
\end{figure}
\end{center}

The present section is devoted to some remarks regarding complementary results and generalizations of the statements of Sections \ref{sezz3} and \ref{sez5}. For simplicity, we will assume $F(1) > 0$ (the case $F(1)=0$ is immediate, as in the previous sections). 
\begin{remark}
\textnormal{
We first give a counterpart of the results of the previous section concerning the behavior of the critical fronts for the equation 
$$
\frac{1}{\eps} \left(\frac{v_\eps'}{\sqrt{1-(v_\eps')^2}}\right)' - c_\eps^* v_\eps' + f(v) =0
$$
($\eps \to 0^+$). Here the picture is much simpler than the one in Section \ref{sez5}: on the one hand, from \eqref{Stimaa} and \eqref{Stimac} it follows immediately that $c_\eps^* \to +\infty$. On the other hand, using the boundedness of $y_\eps$ provided by \eqref{y1}, 
\eqref{costruzionev} immediately implies that $v_\eps'(z) \to 0$ for every $z$. The sequence $v_\eps$ thus converges to the constant $\mathcal{V}_0$.  }
\end{remark}
We now turn to the general equation \eqref{TW2}. In Sections \ref{sezz3} and \ref{sez5}, the use of the parameters $\eps$ and $\gamma$ has allowed us to simplify the statements and to explicitly focus on problems which are extensively studied in literature. However, one could also maintain the two parameters $a$ and $b$ distinct and study the equations (recall \eqref{TW2} and \eqref{notazioney})
\begin{equation}\label{general+}
\frac{1}{a} \left(\frac{v_{a, b}'}{\sqrt{1-\Big(\displaystyle \frac{b}{a}\Big)^2(v_{a, b}')^2}}\right)' - c_{a, b}^* v_{a, b}' + f(v_{a, b}) = 0, \qquad 
\left\{
\begin{array}{l}
\displaystyle y_{a, b}'= c_{a, b}^* \displaystyle \frac{\sqrt{y_{a, b}\left(\displaystyle \frac{2}{a}+\Big(\displaystyle \frac{b}{a}\Big)^2y_{a, b}\right)}}{\displaystyle \frac{1}{a}+\Big(\displaystyle\frac{b}{a}\Big)^2y_{a, b}} - f(v) \vspace{0.2cm}\\
y_{a, b}(0)=0, \, y_{a, b}(1)= 0, \, y_{a, b} > 0 \text{ on } (0, 1).
\end{array}
\right.
\end{equation}
It appears clear, as already suggested by \eqref{Stimaa} and \eqref{Stimac}, that the two key quantities influencing the asymptotic behavior of the critical speed $c_{a, b}^*$ and of the associated critical profile $v_{a, b}$ are $1/a$ and $b/a$; the significant cases are here the ones where at least one between $b/a$ and $1/a$ converges either to $0$ or to $+\infty$. 
An explicit use of \eqref{Stimaa} and \eqref{Stimac}, together with arguments similar to the ones in the previous sections, may then lead to general statements involving both the parameters $a, b$. After having observed that $y_{a, b}$ is uniformly bounded, in the same way as for \eqref{y1}, we just make the following considerations, referring to Figure \ref{comportamenti} (for the reader's convenience, in parentheses we here indicate the type of line used in the figure to give an idea of each situation, keeping in mind that the provided graphic representation is intuitive and only for an illustrative purpose): 
\begin{itemize}
\item (regular black line) if $b/a \to +\infty$, then necessarily one has $c_{a, b}^* \to +\infty$, as a consequence of \eqref{Stimaa} and \eqref{Stimac}; the fact that the solutions of the first differential equation in \eqref{general+} satisfy $\Vert v'_{a, b} \Vert_{L^\infty(\mathbb{R})} \leq a/b \to 0$ then implies, as for Theorems \ref{teorema11eps} and \ref{teoremabis}, that $v_{a, b} \to \mathcal{V}_0$; 
\item (thick black line) if $a \to 0^+$,  $b \to 0^+$ but $a/b \not\to 0$, it is straightforward to check that \eqref{Stimaa} and \eqref{Stimac} produce anyway $c_{a, b}^* \to +\infty$, with order $1/\sqrt{a}$. 
From the second equality in \eqref{costruzionev}, one then deduces that $v_{a, b}' \to 0$ and $v_{a, b} \to \mathcal{V}_0$, using the uniform boundedness of $y_{a, b}$; 
\item (gray line) if $b \to 0^+$ and $a$ remains bounded and bounded away from $0$, the situation can be practically led back to Theorem \ref{modellonuovo}, taking into account Remark \ref{conclusivo}, as well; the limit picture is thus the one of the linear reaction-diffusion equation and $v_{a, b}$ converges in $C^2$ to a front profile $v_{\text{lin}}$ for equation \eqref{lineare}. However, in principle it is not immediate to prove that this front is the critical one, since the limit value $c_{\text{lin}}$ of $c_{a, b}^*$ (which, being bounded, converges up to subsequences) is 
only subject to the bounds \eqref{Stimaa} and \eqref{Stimac};
\item (thick gray line) if $a \to +\infty$ and $b^2/a \to 0$, Proposition \ref{proponuova} yields $c_{a, b}^* \to 0$. As for the critical profiles, the first differential equation in \eqref{costruzionev} and the fact that $y_{a, b}$ is bounded produce $b^2(v'_{a, b})^2/a^2 \to 0$. If $b \to 0^+$, 
since $\Vert y_{a, b} \Vert_{L^\infty(0, 1)} \leq M$ for a suitable constant $M > 0$ the differential equation in \eqref{costruzionev} ensures that 
$
(v_{a, b}^*)'(z) \leq \dfrac{a\sqrt{M(2a+b^2 M)}}{a+b^2 M} = \dfrac{\sqrt{M(2a+b^2 M)}}{1+(b^2/a) M};
$
hence, whether unbounded, $(v_{a, b}^*)'(z)$ diverges at most with order $\sqrt{a}$. 
The first differential equation in \eqref{general+} then provides 
$$
-\frac{1}{\sqrt{a}} \int_{\mathbb{R}} \frac{v_{a, b}'(z)/\sqrt{a}}{\sqrt{1-(b/a)^2(v_{a, b}'(z))^2}} \psi'(z) \, dz + c_{a, b}^* \int_{\mathbb{R}} v_{a, b}(z) \psi'(z) \, dz + \int_{\mathbb{R}} f(v_{a, b}(z)) \psi(z) \, dz = 0
$$
and the fact that $1/a \to 0$ allows one to infer that the first summand vanishes in the limit. Consequently, the limit profile $\bar{v}$ satisfies $\int_{\mathbb{R}} f(\bar{v}(z)) \psi(z) \, dz = 0$. The same conclusion would actually be obtained if $b \not\to0^+$, writing the first differential equation in \eqref{general+} in the form 
$$
-\frac{1}{b} \int_{\mathbb{R}} \frac{\frac{b}{a} v_{a, b}'(z)}{\sqrt{1-(b/a)^2(v_{a, b}'(z))^2}} \psi'(z) \, dz + c_{a, b}^* \int_{\mathbb{R}} v_{a, b}(z) \psi'(z) \, dz + \int_{\mathbb{R}} f(v_{a, b}(z)) \psi(z) \, dz = 0.
$$
In any case, the limit profile $\bar{v}$ coincides with a stepwise function $\mathcal{H}$, but in order to prove that $\mathcal{H} \equiv H$, as in Theorem \ref{conv2}, one should still exclude the convergence to the equilibrium $\alpha$ along an interval (we do not deepen this issue for briefness); 
\item (thickest black line) if $a \to +\infty$, $b^2/a \to +\infty$ and $b/a \not \to +\infty$, one can show that the conclusion of Theorem \ref{ilteorema} holds in the same way. Indeed, one preliminarily notices that the piecewise linear limit profile $\mathcal{V}_L$, whenever nonzero, will have slope equal to $a/b$. 
Repeating the argument in the proof of Proposition \ref{pezzolineare}, if the limit profile contains a piecewise linear part with slope $a/b$ whenever taking values between $v_0$ and $v_1$, one then finds that 
$$
\bar{c}=\frac{b}{a}\left(\frac{F(v_1)-F(v_0)}{v_1-v_0} +\frac{\bar{y}(v_1)-\bar{y}(v_0)}{v_1-v_0}\right). 
$$
Similarly as for the proof of Lemma \ref{lemma1b}, taking into account \eqref{Stimaa} it follows that the critical fronts $v_{a, b}$ cannot converge to the fully piecewise linear profile $\mathcal{V}_L$, while the regular scenario is excluded proceeding similarly as in Lemma \ref{lemma1}. The picture is thus analogous to the one for \eqref{TWdelta}.
\end{itemize}
Briefly speaking, one may then expect to find all the records encountered in the previous sections, according to the trend of the parameters $a$ and $b$. The only case which seems subtler to deal with and may deserve further attention is $a \to +\infty$, $b \to +\infty$ with $b/a \to 0$ but $b^2/a \not \to 0$; here $c_{a, b}^* \to 0$ but the shape of the limit profile is less clear from the previous computations (we qualitatively indicate the corresponding situation in Figure \ref{comportamenti} with a dashed gray line). However, the numerical simulation shown in Figure \ref{convergenzalast} seems to suggest that the limit profile is again stepwise (in the example, coinciding with $H$); here, for $b=\sqrt{a}$ (left), $b=a^{3/4}$ (middle), $b=a^{4/5}$ (right) and $a=10^{-3}$ (dashed), $a=10^{-5}$ (black), we display the approximated shape of the critical front, respectively obtained for $c^*_{a, b} \approx 0.06/0.006$ (left), $c^*_{a, b} \approx 0.06/0.011$ (middle), $c^*_{a, b} \approx 0.064/0.019$ (right). Notice that, the more $b$ approaches $a$, the more the profile is still reminiscent of the shape of $\mathcal{V}_{L, I}$ for $a=10^{-5}$, and the convergences $c^*_{a, b} \to 0$ and of $v_{a, b} \to H$ are expected to be slower, reasonably in accord with the observations in Section \ref{sez5}. 

\begin{center}
\begin{figure}[!h]
\makebox[\linewidth][c]{
\begin{tikzpicture}[xscale=0.9,yscale=0.9]

\draw[thick]{(-1, 0)--(12,0)};
\draw[thick]{(0, -1)--(0, 12)};
\node at (12.3, -0.3){$a$};
\node at (-0.3, 12.3){$b$};

\draw[line width=0.05]{(0, 0)--(12, 12)};

\draw[very thick, domain=0:2.5] plot(\x, 0.111*\x*\x) [arrow outside={}{0.5}];
\draw[thick, domain=0:0.75] plot ({\x}, {sqrt(9*\x)}) [arrow outside={}{0.5}];
\draw[very thick, domain=0:2] plot(\x, \x) [arrow outside={}{0.5}];
\node at (2, 3.2){$\begin{array}{l} c_{a, b}^* \to +\infty  \\ v_{a, b} \to \mathcal{V}_0 \end{array} $};

\draw[thick, domain=0:2.5] plot(\x, 5.6) [arrow outside={}{0.5}];
\draw[thick, domain=0:2.5] plot (\x, 0.1*\x+5.8) [arrow outside={}{0.5}];
\draw[thick, domain=0:2.5] plot(\x, -0.1*\x+5.4) [arrow outside={}{0.5}];
\node at (4,5.65){$\begin{array}{l} c_{a, b}^* \to +\infty  \\ v_{a, b} \to \mathcal{V}_0 \end{array} $};

\draw[thick, domain=0.1:2.5] plot(\x, 0.35/\x+8.5) [arrow outside={}{0.4}];
\draw[thick, domain=0.16667:2.5] plot(\x, 0.5/\x+9) [arrow outside={}{0.4}];
\draw[thick, domain=0.2:2.5] plot(\x, 0.4/\x+10) [arrow outside={}{0.4}];
\node at (2.5, 11.5){$\begin{array}{l} c_{a, b}^* \to +\infty  \\ v_{a, b} \to \mathcal{V}_0 \end{array} $};

\draw[gray, domain=0:2.8] plot(6.1, \x) [arrow outside={}{0.5}];
\draw[gray, domain=5.5:5.9] plot(\x, 7*\x-38.5) [arrow outside={}{0.5}];
\draw[gray, domain=6.7:6.3] plot(\x, -7*\x+46.9) [arrow outside={}{0.5}];
\node at (6.1, 3.7){\textcolor{gray}{$\begin{array}{l} c_{a, b}^* \to c_{\text{lin}}  \\ v_{a, b} \to v_{\text{lin}}\end{array} $}};

\draw[thick, domain=12:9.2] plot(6.1, \x) [arrow outside={}{0.5}];
\draw[thick, domain=5.919:5.75] plot(\x, {(\x-5.05)/(6.1-\x)+7.2}) [arrow outside={}{0.5}];
\draw[thick, domain=6.281:6.45] plot(\x, {(7.15-\x)/(\x-6.1)+7.2}) [arrow outside={}{0.5}];
\node at (6.1, 8.5){$\begin{array}{l} c_{a, b}^* \to +\infty  \\ v_{a, b} \to \mathcal{V}_0 \end{array} $};

\draw[gray, very thick, domain=12:9.5] plot(\x, 5.6) [arrow outside={}{0.5}];
\draw[gray, very thick, domain=12:9.5] plot (\x, {0.01*(6*\x+2)/(\x-8)+5.7}) [arrow outside={}{0.5}];
\draw[gray, very thick, domain=12:9.5] plot (\x, {0.01*(6*\x+2)/(8-\x)+5.5}) [arrow outside={}{0.5}];
\node at (8, 5.65){\textcolor{gray}{$\begin{array}{l} c_{a, b}^* \to 0  \\ v_{a, b} \to \mathcal{H} \end{array}$}};

\draw[thick, domain=10.41:9.5] plot(\x, {(\x-9)^2+10}) [arrow outside={}{0.5}];
\draw[line width=1.5, domain=12:9.5] plot (\x, \x) [arrow outside={}{0.5}];
\draw[gray, very thick, domain=12:9.5] plot(\x, {(\x-9)^(1/4)+7.8}) [arrow outside={}{0.5}];
\node at (9.2, 12.5){$\begin{array}{l} c_{a, b}^* \to +\infty  \\ v_{a, b} \to \mathcal{V}_0 \end{array} $};
\node at (12.5, 8.3){\textcolor{gray}{$\begin{array}{l} c_{a, b}^* \to 0  \\ v_{a, b} \to \mathcal{H} \end{array} $}};

\draw[gray, dashed, domain=12:9.5] plot(\x, {sqrt(\x-9)+8.25}) [arrow outside={}{0.5}];
\node at (13.1, 9.8){\textcolor{gray}{$c_{a, b}^* \to 0$}};

\node at (13.85, 12.5){{$c_{a, b}^* \to \bar{c} > 0$}};
\node at (15, 11.95){{$v_{a, b} \to \bar{v}$ \text{ as in Section \ref{sez5}}}};
\node at (14.65, 11.4) {{($\bar{v}=\mathcal{V}_L$ or $\bar{v}=\mathcal{V}_{L, I}$)}};

\draw[gray, very thick, domain=12:9.5] plot(\x, {1/(\x-6)}) [arrow outside={}{0.5}];
\draw[gray, very thick, domain=12:9.5] plot(\x, {2/(\x-7.4)}) [arrow outside={}{0.5}];
\draw[gray, very thick, domain=12:9.5] plot(\x, {3/(\x-8)}) [arrow outside={}{0.5}];
\node at (10, 3.2){\textcolor{gray}{$\begin{array}{l} c_{a, b}^* \to 0  \\ v_{a, b} \to \mathcal{H} \end{array} $}};

\end{tikzpicture}}
\caption{\small{A qualitative snapshot of the behavior of the critical front profiles as $a$ and $b$ vary}} 
\label{comportamenti}
\end{figure}
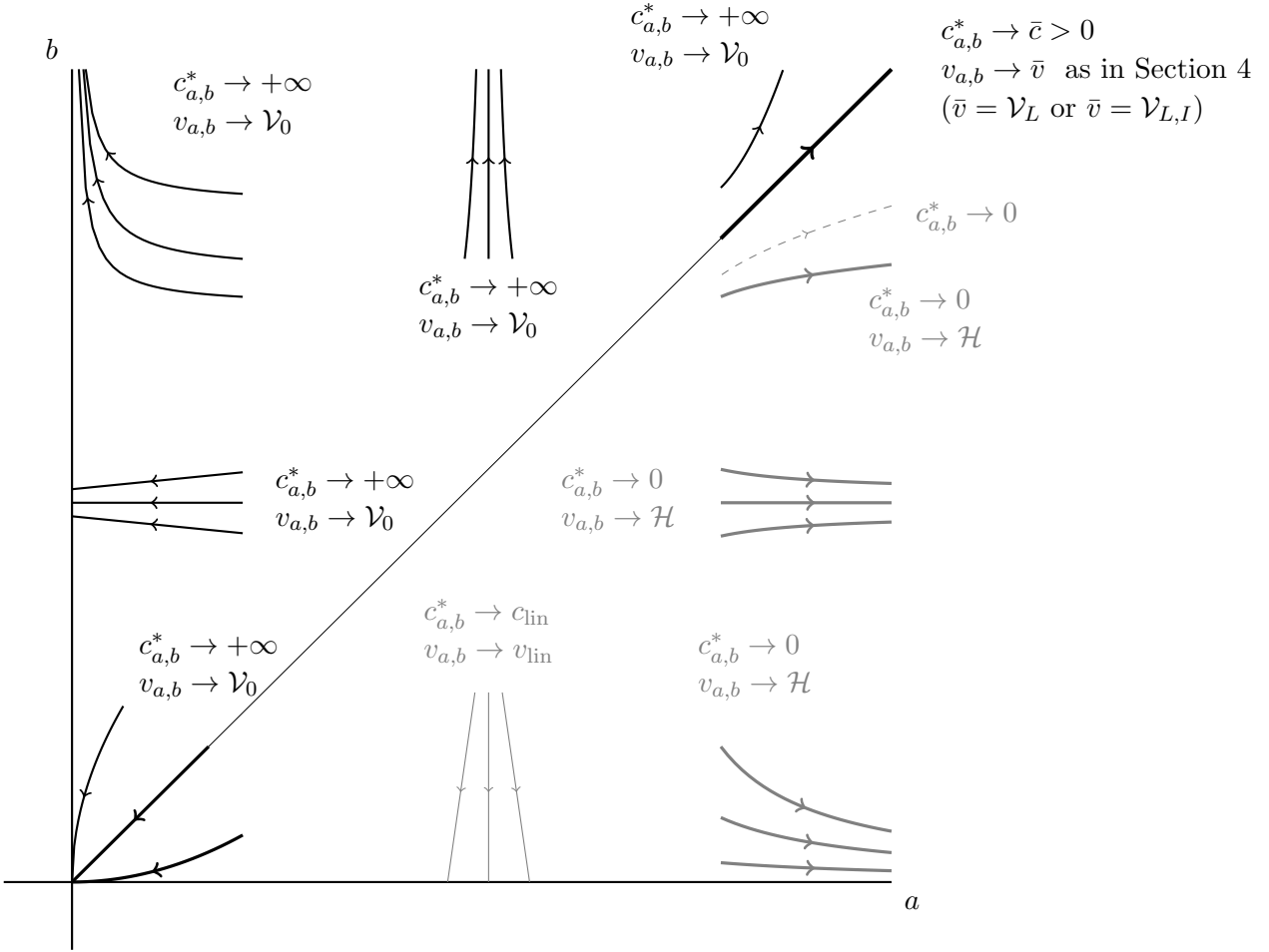
\end{center}

Figure \ref{comportamenti} evidently shows the importance of the interplay between the parameters $a$ and $b$ in the limit outcome. From a visual point of view, the transition between constant and jumping limit profiles can occur in two different ways: when moving increasingly along the $a$-axis, the limit profile passes from being constant to regular (coinciding with the one of the linear diffusion equation) until this regular shape ``breaks'', becoming stepwise. On the other hand, such a transition may be encountered also at infinity ($a, b \to +\infty$), and in this case one-sided sharp (or possibly fully piecewise linear, in the type C case) limit profiles are found (for instance, for $a=b$), then ``breaking'' again to a piecewise constant function. In both cases, the limit speed decreases monotonically from $+\infty$ to $0$, passing through finite values. As a general behavior, we observe moreover that the more the profile becomes singular, the slower becomes the limit speed.

\section{Conclusion}

We have analyzed 
traveling front-type solutions for the $1$-dimensional PDE 
$$
u_t=\left(\frac{u_x}{\sqrt{a^2-b^2 u_x^2}}\right)_x  + f(u),
$$
in dependence on the parameters $a$ and $b$. We have worked on providing readable \textbf{estimates of the critical speed} (complementing existing ones) and on the analysis of the \textbf{limit behavior of the critical front profile}. 
Among the physically predicted results, we have seen that for $b \to 0^+$ and bounded $a$ one has the \textbf{convergence to the critical profile of the linear reaction-diffusion model}, together with the \textbf{convergence of the relative propagation speed}: this confirms that Maxwell's theory is recovered from the Born-Infeld one for maximal field strength going to infinity. On the other hand, $b/a \to +\infty$ always leads to find a \textbf{constant limit profile propagating with infinite speed}, since intuitively $u_x$ has to stay arbitrarily small in order for the second-order operator to be well-defined. Finally, if $a \to +\infty$ sufficiently fast with respect to $b$, then the profiles have arbitrarily large room to increase their steepness and a \textbf{stepwise limit profile} is found.
A new surprising behavior arises instead for the \textbf{singular perturbation problem}, with the appearance of a \textbf{limit front profile sharpening on one side only} (for KPP or combustion reaction terms), a feature which seems new in the theory of traveling fronts in case of nondegenerate diffusivity. Figure \ref{comportamenti} illustrates the discussed behaviors for general $a, b$.

\section*{Appendix: numerical approximations of the critical speed}

We here list the numerical approximations of the critical speed in the above figures. Our shooting procedure has been implemented determining $y_{c, f}^-$ numerically in dependence on $c$, starting from the lower bound in \eqref{Stimaa} and \eqref{Stimac} on, until $y_{c, f}^-(0)$ and $(y_{c, f}^-)'(0)$ were sufficiently close to $0$. Despite our approximations of $c^*$ often contained more than three digits after decimal, we here only report the first three. We remark that these are only approximations of the real critical speed; to find completely reliable values of $c^*$, one should probably proceed through some computer-assisted arguments. 
\begin{table}[ht]
\begin{center}
\begin{tabular}{|c|c|c|c|c|c|c|}
\hline
& $\gamma$ & $1$ & $5$ & $10$ & $50$ & $100$ \\
\hline
Figure \ref{convergenza1} & $c^*_\gamma$ & \text{\textsc{0.142}} & \text{\textsc{0.167}} & \text{\textsc{0.227}} & \text{\textsc{0.874}} &  \text{\textsc{1.710}}  \\
\hline 
\hline
& $\gamma$ & $10^{-4}$ & $10^{-2}$ & $1$ & $10$ & $100$ \\
\hline
Figure \ref{convergenza2} & $c^*_\gamma$ & \text{\textsc{0.020}} & \text{\textsc{0.201}}  & \text{\textsc{1.893}}  & \text{\textsc{5.919}}  & \text{\textsc{22.901}}  \\
\hline 
\hline
& $\eps^2$ & $10^{-1}$ & $10^{-2}$ & $10^{-3}$ & $10^{-4}$ & $10^{-6}$ \\
\hline
Figure \ref{fisher} & $c^*_\eps$ & \text{\textsc{1.011}} & \text{\textsc{0.569}} & \text{\textsc{0.332}} & \text{\textsc{0.229}} & \text{\textsc{0.193}}   \\
\hline
Figure \ref{nagy} & $c^*_\eps$ & \text{\textsc{1.318}} & \text{\textsc{0.875}} & \text{\textsc{0.707}} & \text{\textsc{0.652}} & \text{\textsc{0.626}}  \\
\hline
Figures \ref{tipoB1}-\ref{convergenzacsm} & $c^*_\eps$ & \text{\textsc{0.284}} & \text{\textsc{0.167}} & \text{\textsc{0.106}} & \text{\textsc{0.078}} & \text{\textsc{0.063}} \\
\hline
Figure \ref{convergenzacpl} & $c^*_\eps$ & \text{\textsc{0.040}} & \text{\textsc{0.024}} & \text{\textsc{0.015}} & \text{\textsc{0.011}} & \text{\textsc{0.009}} \\
\hline
\end{tabular}
\caption{Approximations of the critical speed in the indicated figures}
\end{center}
\end{table}


{\footnotesize
\begin{thebibliography}{99}
\bibitem{AroWei}{D. G. Aronson and H. F. Weinberger, \emph{Multidimensional nonlinear diffusion arising in population genetics}, Adv. Math.
\textbf{30} (1976), 33--76.}
\bibitem{Azz06}{A. Azzollini, \emph{On a prescribed mean curvature equation in Lorentz-Minkowski space}, J. Math. Pures Appl. (9) \textbf{106} (2016), 1122--1140.}
\bibitem{BarSim}{R. Bartnik and L. Simon, \emph{Spacelike hypersurfaces with prescribed boundary values and mean curvature},
Comm. Math. Phys. \textbf{87} (1982/83), 131--152.}
\bibitem{BerNir}{H. Berestycki and L. Nirenberg, \emph{Travelling fronts in cylinders}, Ann. Henri Poincar\'e \textbf{9} (1992), 497--572.}
\bibitem{BonCoeNys}{D. Bonheure, I. Coelho and M. Nys, \emph{Heteroclinic solutions of singular quasilinear bistable equations}, NoDEA Nonlinear Differential Equations Appl. \textbf{24} (2017), Paper No. 2, 29 pp.}
\bibitem{BondAvPomRei}{D. Bonheure, P. d'Avenia, A. Pomponio and W. Reichel, \emph{Equilibrium measures and equilibrium potentials in the Born-Infeld model}, J. Math. Pures. Appl. \textbf{139} (2020), 35--62.}
\bibitem{BonSan}{D. Bonheure and L. Sanchez, {\it Heteroclinic orbits for some classes of second and fourth order differential equations}, 
In:
A. Canada, P. Dr\'abek, A. Fonda eds., Handbook of Differential Equations: Ordinary Differential
Equations, vol. 3,
Elsevier, Amsterdam, 2006, 103--202.}
\bibitem{BorInf}{M. Born and L. Infeld, \emph{Foundations of the new field theory}, Proc. R. Soc. Lond. Ser. A \textbf{144} (1934), 425--451.}
\bibitem{BosColNor}{A. Boscaggin, F. Colasuonno and B. Noris, \emph{Positive radial solutions for the Minkowski-curvature equation with Neumann boundary conditions}, Discrete Contin. Dyn. Syst. Ser. S \textbf{13} (2020), 1921--1933.}
\bibitem{BosFel}{A. Boscaggin and G. Feltrin, \emph{Positive periodic solutions to an indefinite
Minkowski-curvature equation}, J. Differential Equations \textbf{269} (2020), 5595--5645.}
\bibitem{CarKie}{H. Carley and K.-H. Kiessling, \emph{Constructing graphs over $\mathbb{R}^n$ with small prescribed mean-curvature}, Math. Phys. Anal. Geom. \textbf{18} (2015), 25 pp.}
\bibitem{CoeSan}{I. Coelho and L. Sanchez, \emph{Travelling wave profiles in some models with nonlinear diffusion}, Appl. Math. Comp. \textbf{235} (2014), 469--481.} 
\bibitem{CorObeOmaRiv}{C. Corsato, F. Obersnel, P. Omari and S. Rivetti, \emph{Positive solutions of the Dirichlet problem
for the prescribed mean curvature equation in Minkowski space}, J. Math. Anal. Appl. \textbf{405}
(2013), 227--239.}
\bibitem{Diek}{O. Diekmann, \emph{Limiting behaviour in an epidemic model}, Nonlinear Anal. \textbf{1} (1976/77), 459--470.}
\bibitem{DraTak}{P. Dr\'abek and P. Taka\v{c}, \emph{Travelling waves in the Fisher-KPP equation with nonlinear diffusion and a non-Lipschitzian reaction term}, preprint, ArXiv 1803.10306}
\bibitem{FifMcL}{P. C. Fife and J. B. Mc Leod, \emph{The approach of solutions of nonlinear diffusion equations to travelling front
solutions}, Arch.
Ration. Mech. Anal. \textbf{65} (1977), 335--361.}
\bibitem{FolPlaStr}{R. Folino, R. Plaza and M. Strani, \emph{Metastable patterns for a reaction-diffusion model with mean curvature-type diffusion}, J. Math. Anal. Appl. \textbf{493} (2021), 124455.}
\bibitem{Gar}{M. Garrione, \emph{Vanishing diffusion limits for planar fronts in bistable models with saturation}, Trans. Amer. Math. Soc. \textbf{374} (2021), 3999-4021.}
\bibitem{GarSan}{M. Garrione and L. Sanchez, \emph{Monotone traveling waves for reaction-diffusion equations involving the curvature operator}, Bound. Value Probl. \textbf{2015:45}, 1--31.}
\bibitem{Gib}{G.W. Gibbons, \emph{Aspects of Born-Infeld theory and string/M-theory}, Rev. Mex. F\'is. S  \textbf{1} (2003), 19--29.}
\bibitem{GilKer}{B.H. Gilding and R. Kersner, Travelling waves in Nonlinear Diffusion-Convection Reaction, Birkh\"auser, Basel, 2004.}
\bibitem{HadRot}{K.P. Hadeler and F. Rothe, \emph{Travelling fronts in nonlinear diffusion equations}, J. Math. Biol. \textbf{2} (1975), 251--263.}
\bibitem{HilKim}{D. Hilhorst and Y-J. Kim, \emph{Diffusive and inviscid traveling waves of the Fisher equation and nonuniqueness of wave speed}, Appl. Math. Lett. \textbf{60} (2016), 28--35.}
\bibitem{HodHux}{A.L. Hodgkin and A.F. Huxley, \emph{A quantitative description of membrane current and its application to conduction and excitation in nerve}, J. Physiol. \textbf{117} (1952), 500--544.}
\bibitem{Kru10}{S.I. Kruglov, \emph{On generalized Born-Infeld electrodynamics}, J. Phys. A: Math. Theor. \textbf{43} (2010), 375402, 1--8.}
\bibitem{Kru2}{S.I. Kruglov, \emph{Notes on Born-Infeld-type electrodynamics}, Mod. Phys. Lett. A \textbf{32} (2017), No. 36, 1750201.}
\bibitem{MalMar}{L. Malaguti and C. Marcelli, \emph{Travelling wavefronts in reaction-diffusion equations with convection effects and non-regular
terms}, Math. Nachr. \textbf{242} (2002), 148--164.}
\bibitem{MalMar03}{L. Malaguti and C. Marcelli, \emph{Sharp profiles in degenerate and doubly degenerate Fisher-KPP equations}, J. Differential Equations \textbf{195} (2003), 471--496.}
\bibitem{Maw}{J. Mawhin, \emph{Nonlinear boundary value problems involving the extrinsic mean curvature operator}, 
Math. Bohem. \textbf{139} (2014), 299--313.} 
\bibitem{Nag75}{T. Nagylaki, \emph{Conditions for the existence of clines}, Genetics \textbf{3} (1975), 595--615.}
\bibitem{RacStaTvr06}{I. Rach{\r u}nkov\'a, S. Stan\v{e}k and M. Tvrd\'y, Singularities and Laplacians in boundary value problems for nonlinear ordinary differential equations, Handbook of differential equations: ordinary differential equations. Vol. III, 607–722, Elsevier, Amsterdam, 2006.}
\bibitem{Ros}{P. Rosenau, \emph{Free-energy functionals at the high-gradient limit}, Phys. Rev. A
\textbf{41} (1990),
2227--2230.
}
\bibitem{SanMai94}{F. S\'anchez-Gardu\~no and P.K. Maini, \emph{Existence and uniqueness of a sharp travelling wave in degenerate non-linear diffusion Fisher-KPP equations}, J. Math. Biol. \textbf{33} (1994), 163--192.}
\bibitem{Wal}{W. Walter, Ordinary Differential Equations, Springer, New York, 1998.}
\bibitem{Yan}{Y. Yang, \emph{Classical solutions in the Born-Infeld theory}, Proc. Roy. Soc. London Ser. A \textbf{456} (2000), 615--640.}

\end{thebibliography}
}
\end{document}